\documentclass[11]{amsart}
\usepackage[small]{caption}

\usepackage{amsmath,amsfonts,amscd,amssymb,epsf, epsfig}\usepackage{graphicx}

\begin{document}

\title[Generalized Whittle-Mat$\acute{\text{E}}$rn random field]{  Generalized Whittle-Mat$\acute{\text{E}}$rn random field as a model of correlated fluctuations }

\author{S.~C.~Lim}\email{sclim@mmu.edu.my}
\address{Faculty of Engineering,
Multimedia University, Jalan Multimedia, Cyberjaya, 63100, Selangor
Darul Ehsan, Malaysia.}

\author{L.~P.~Teo}\email{lpteo@mmu.edu.my}
\address{Faculty of Information
Technology, Multimedia University, Jalan Multimedia, Cyberjaya,
63100, Selangor Darul Ehsan, Malaysia.}

\begin{abstract}
This paper considers a generalization of Gaussian random field with covariance function of Whittle-Mat$\acute{\text{e}}$rn family. Such a random field can be obtained as the solution to the fractional stochastic differential equation with two fractional orders.   Asymptotic properties of the covariance functions belonging to this generalized Whittle-Mat$\acute{\text{e}}$rn family are studied, which are used to deduce the sample path properties of the random field. The Whittle-Mat$\acute{\text{e}}$rn field has been widely used in modeling geostatistical data such as sea beam data, wind speed, field temperature and soil data. In this article we show that generalized Whittle-Mat$\acute{\text{e}}$rn field provides a more flexible model for wind speed data.
\end{abstract}


\maketitle

\section{Introduction}
Random fields play an important role in geostatistics, which deals with problems stretching from the resource evaluation such as the estimation of ore resources in mining and oil deposits in oil exploration, pollution evaluation in environmental sciences, to hydrology, meteorology, agriculture, etc. \cite{1,2,3}. For examples, environmental resource models carry out spatial statistical analysis in the quantity of resources available such as the volume of available water, forest, etc., or their quality such as concentration of contaminants in air, water or soil samples. Random fields and their covariance functions or equivalently their variograms are used widely in the modeling of observed spatial data as these data are likely to be spatially dependent. The earlier developments of the subject include work by Whittle \cite{4,5}, Mat$\acute{\text{e}}$rn \cite{6,7}, Tatarski \cite{8}, Matheron \cite{9} and others. The Gaussian random fields defined using the covariance functions from the Whittle-Mat$\acute{\text{e}}$rn (WM) covariance class are widely used to model isotropic spatial processes in two and three dimensions.

The WM class of covariance functions \cite{10, 11, 12, 13, 41} has recently received considerable interest in geostatistics due to its great flexibility for modeling the spatial variations, in particular its ability to model behaviors of empirical variogram near the origin. Unlike other popular covariance models, the WM model has a parameter that characterizes the smoothness of the associated random field. Due to this reason, Stein \cite{14} strongly recommended the WM class for the modeling of spatial covariance.

A special case of the WM model was first obtained by Whittle \cite{4}, who showed that a Gaussian random field with covariance function belonging to WM class can be obtained as a solution to a stochastic differential equation. The general form of WM model was given by Mat$\acute{\text{e}}$rn \cite{6} and Tatarski \cite{8}, and was also considered by Matheron \cite{9} and Shkarofsky \cite{15}. It can be associated with von K$\acute{\text{a}}$rm$\acute{\text{a}}$n spectrum \cite{16, 17} in the modeling of wind speed. A comprehensive historical account on the WM class was given by Guttorp and Gneiting \cite{18}, who first called such a class of covariance functions as WM covariance family, but they later changed it to Mat$\acute{\text{e}}$rn covariance family \cite{19}.

Recall that the smoothness of a random field is characterized by the fractal dimension, a local property which is determined by the asymptotic properties of the covariance near zero lag and its value depends on the smoothness parameter of the WM covariance. On the other hand, the strength of the spatial correlation is determined by the scale parameter  and for large time lag it decays exponentially. In this paper, we proposed a new generalization of WM covariance class with an additional parameter which plays the role of scale or memory parameter, and the spatial correlation strength for large time lag now varies hyperbolically, with exponential decay as a special case.

 In the next section we recall some basic facts on the WM covariance class and the random field associated with it. This will be followed by the introduction of the GWM covariance class and the corresponding random field. The asymptotic properties of the GWM covariance function are studied in section \ref{sec3}. Based on these properties we are able to obtain the fractal dimension of the graph of the random field in GWM model. This random field satisfies a weaker self-similar property called local self-similarity, and it is short-range dependent. Simulations of the GWM covariance function and the random field (in two dimensions) are given. In the subsequent section, GWM process is applied to model  wind speed and compared to model provided by WM process. Other possible generalizations and applications of the GWM random field are discussed in the concluding section.

\section{Generalized Whittle-Mat$\acute{\text{e}}$rn Model}\label{sec2}
The WM class of covariance functions is given by \cite{1,2,3, 14}:
\begin{equation}\label{eq1}
C(\boldsymbol{t})=\frac{2^{1-\frac{n}{2}-\gamma}}{\pi^{\frac{n}{2}}\Gamma(\gamma)}
\left(\frac{|\boldsymbol{t}|}{\lambda}\right)^{\gamma-\frac{n}{2}}K_{\gamma-\frac{n}{2}}(\lambda|\boldsymbol{t}|),
\end{equation}where $K_{\nu}(z)$   is the modified Bessel function of second kind (or Macdonald function), $\boldsymbol{t}\in\mathbb{R}^n$, $|\boldsymbol{t}|=\sqrt{t_1^2+\ldots+t_n^2}$ is the Euclidean norm of $\boldsymbol{t}$, $\lambda>0$   is a scale parameter controlling the spatial range of the covariance, and   $\nu=\gamma-(n/2)>0$ is the smoothness parameter governing the level of smoothness of the associated Gaussian random field $Y(\boldsymbol{t})$. Note that the WM covariance \eqref{eq1} has the same functional form as the characteristic function of the multivariate $t$--distribution \cite{20}. The spectral density of   $Y(\boldsymbol{t})$ is given by the Fourier transform of \eqref{eq1}:
\begin{equation}\label{eq2}\begin{split}
S(\boldsymbol{\omega})=F(C(\boldsymbol{t}))=&\frac{1}{(2\pi)^n}\int_{\mathbb{R}^n} C(\boldsymbol{t}) e^{-i\boldsymbol{\omega}.\boldsymbol{t}}d^n\boldsymbol{t}
=\frac{1}{(2\pi)^n}\frac{1}{\left(|\boldsymbol{\omega}|^2+\lambda^2\right)^{\gamma}}.\nonumber
\end{split}\end{equation}The Gaussian random field $Y(\boldsymbol{t})$  with covariance \eqref{eq1} can be obtained as the solution to the following fractional stochastic differential equation \cite{4}:
\begin{equation}\label{eq3}
\left(-\Delta +\lambda^2\right)^{\frac{\gamma}{2}}Y(\boldsymbol{t})=\eta(\boldsymbol{t}),
\end{equation}where $\Delta =\frac{\partial^2}{\partial t_1^2}+\ldots+\frac{\partial^2}{\partial t_n^2}$ is the $n$-dimensional Laplacian, and $\eta(\boldsymbol{t})$  is the standard white noise defined by
\begin{equation}\label{eq4}
\langle \eta(\boldsymbol{t})\rangle =0,\hspace{1cm}\langle \eta(\boldsymbol{t})\eta(\boldsymbol{s})\rangle =\delta(\boldsymbol{t}-\boldsymbol{s}).
\end{equation}
One early generalization of WM family of covariance functions was proposed by Shkarofsky \cite{15}. Based on the argument that a covariance function for turbulence needs to have no cusp, it is required to have zero derivative at the origin and a second derivative that is finite and negative. In order to satisfy these requirements, he generalized \eqref{eq1} to a covariance with two complementary parameters:
\begin{equation}\label{eq5}
C(\boldsymbol{t})=\frac{\left(\lambda\sqrt{|\boldsymbol{t}|^2+\xi^2}\right)^{\nu}K_{\nu}\left(\lambda\sqrt{|\boldsymbol{t}|^2+\xi^2}\right)}{(\lambda\xi)^{\nu}K_{\nu}(\lambda \xi)}.\end{equation}
 Clearly, up to constants, \eqref{eq5} reduces to the covariance in WM class \eqref{eq1}  when $\xi\rightarrow 0^+$.  There also exist generalizations of WM class to a non-stationary class of covariance functions that allows for anisotropy, one such generalization is \cite{21}:
\begin{equation}
C(t_1, t_2)=\left(\frac{\lambda(t_1+t_2)}{2}\right)^{-\nu}K_{\nu}\left(2\sqrt{\frac{\lambda(t_1+t_2)}{2}}\right).
\end{equation}

	In view of the wide applications of fractal operators in physics \cite{54}, we propose another generalization of WM class of covariance function by extending the fractional stochastic differential equation \eqref{eq3} to one with two fractional orders:
\begin{equation}\label{eq7}
\left[\left(-\Delta\right)^{\alpha}+\lambda^2\right]^{\frac{\gamma}{2}}Y_{\alpha,\gamma}(\boldsymbol{t})=\eta(\boldsymbol{t}),
\end{equation}
with $\lambda, \gamma>0$   and $\alpha\in (0, 1]$, and the Riesz fractional derivative $\mathbf{D}^{2\alpha}=(-\Delta)^{\alpha}$  is defined by:
\begin{equation}\label{eq8}\mathbf{D}^{2\alpha}f=(-\Delta)^{\alpha}f=F^{-1}\left\{|\boldsymbol{\omega}|^{2\alpha}F[f](\boldsymbol{\omega})\right\}\end{equation}	    or    \begin{equation} (F\mathbf{D}^{2\alpha}f)(\boldsymbol{\omega})=|\boldsymbol{\omega}|^{2\alpha}(F[f])(\boldsymbol{\omega}),\end{equation}	
 One can regard the fractional operator $\left[\left(-\Delta\right)^{\alpha}+\lambda^2\right]^{\frac{\gamma}{2}}$  as a "shifted" Riesz derivative and it has formally   the series representation:
\begin{equation}
\left[\left(-\Delta\right)^{\alpha}+\lambda^2\right]^{\frac{\gamma}{2}}=\sum_{j=1}^{\infty} \begin{pmatrix} \gamma/2\\j\end{pmatrix}\lambda^{\gamma-2j}(-\Delta)^{\alpha j}.
\end{equation}
See reference \cite{24} for a more rigorous treatment of this operator based on hypersingular integrals.  Now by using
\begin{equation}\left(\left[\left(-\Delta\right)^{\alpha}+\lambda^2\right]^{\frac{\gamma}{2}}f\right)(\boldsymbol{t})=
F^{-1} \left( \left[|\boldsymbol{\omega}|^{2\alpha}+
\lambda^2\right]^{\frac{\gamma}{2}}F[f](\boldsymbol{\omega})\right)(\boldsymbol{t}),\end{equation}
the solution to \eqref{eq7} is found to be
\begin{equation}\label{eq11}
Y_{\alpha,\gamma}(\boldsymbol{t})=\frac{1}{(2\pi)^{\frac{n}{2}}}\int\limits_{\mathbb{R}^n}
\frac{e^{it.\boldsymbol{\omega}}\hat{\eta}(\boldsymbol{\omega})}{\left(|\boldsymbol{\omega}|^{2\alpha}+\lambda^2\right)^{\frac{\gamma}{2}}}
d^n\boldsymbol{\omega},
\end{equation}where $\hat{\eta}(\boldsymbol{\omega})=F[\eta](\boldsymbol{\omega})$ is the Fourier transform of the white noise. For convenience, we call  $Y_{\alpha,\gamma}(\boldsymbol{t})$ the GWM  (generalized Whittle-Mat$\acute{\text{e}}$rn)  field.
The representation \eqref{eq11} shows that the GWM field $Y_{\alpha,\gamma}(\boldsymbol{t})$ is a centered Gaussian field with covariance function\begin{figure}  \epsfxsize=0.49\linewidth
\epsffile{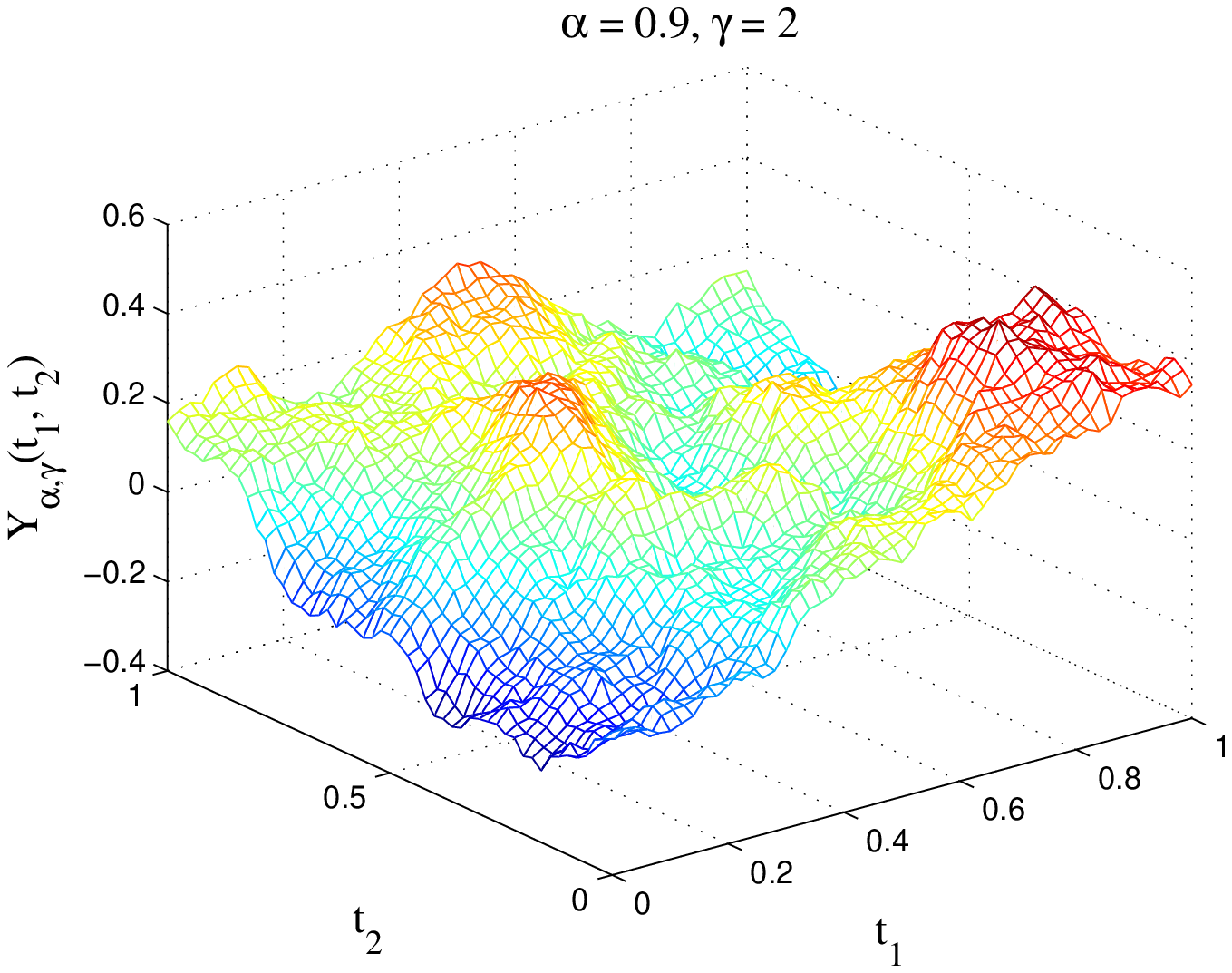} \epsfxsize=0.49\linewidth
\epsffile{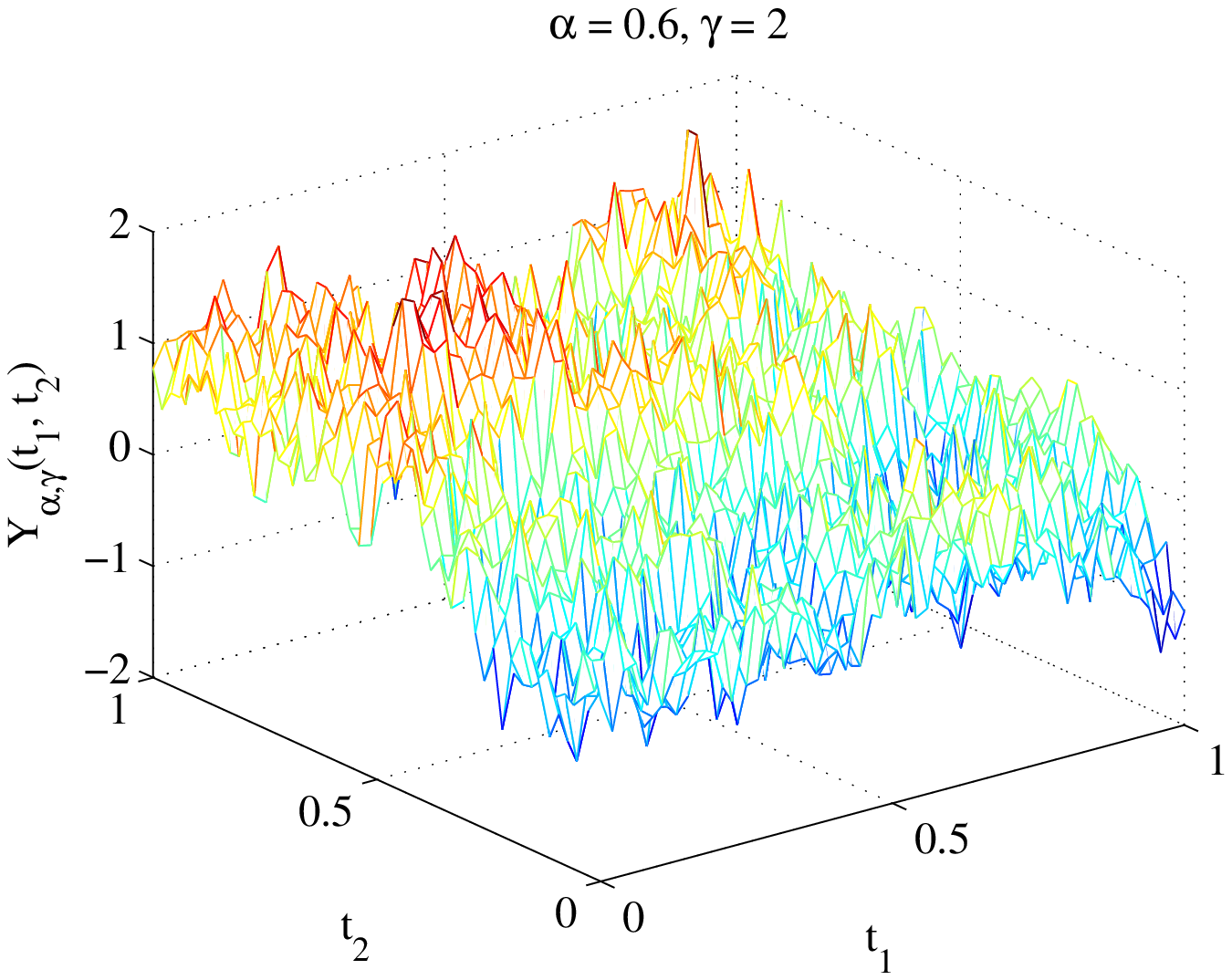}\\
  \epsfxsize=0.49\linewidth
\epsffile{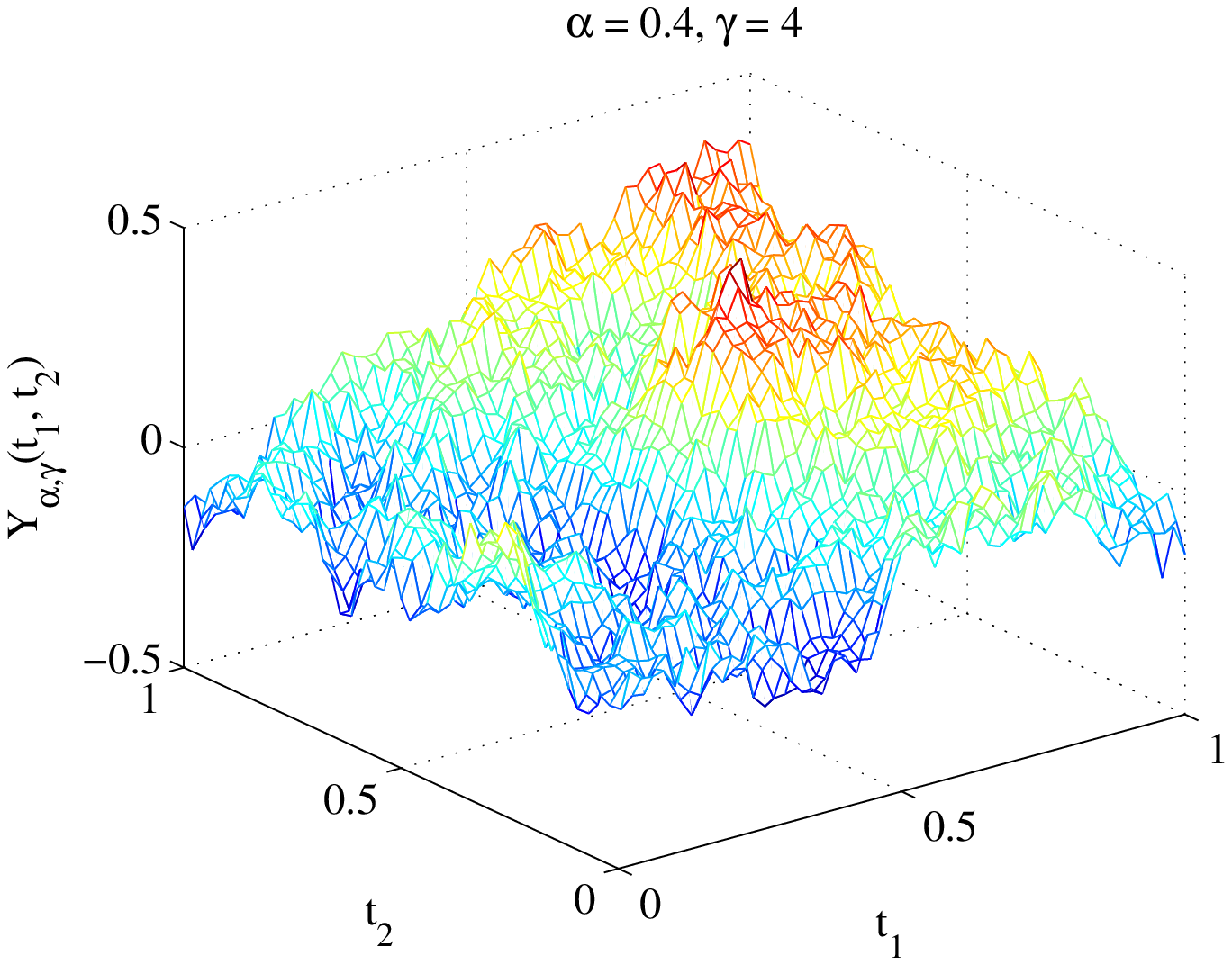}  \epsfxsize=0.49\linewidth
\epsffile{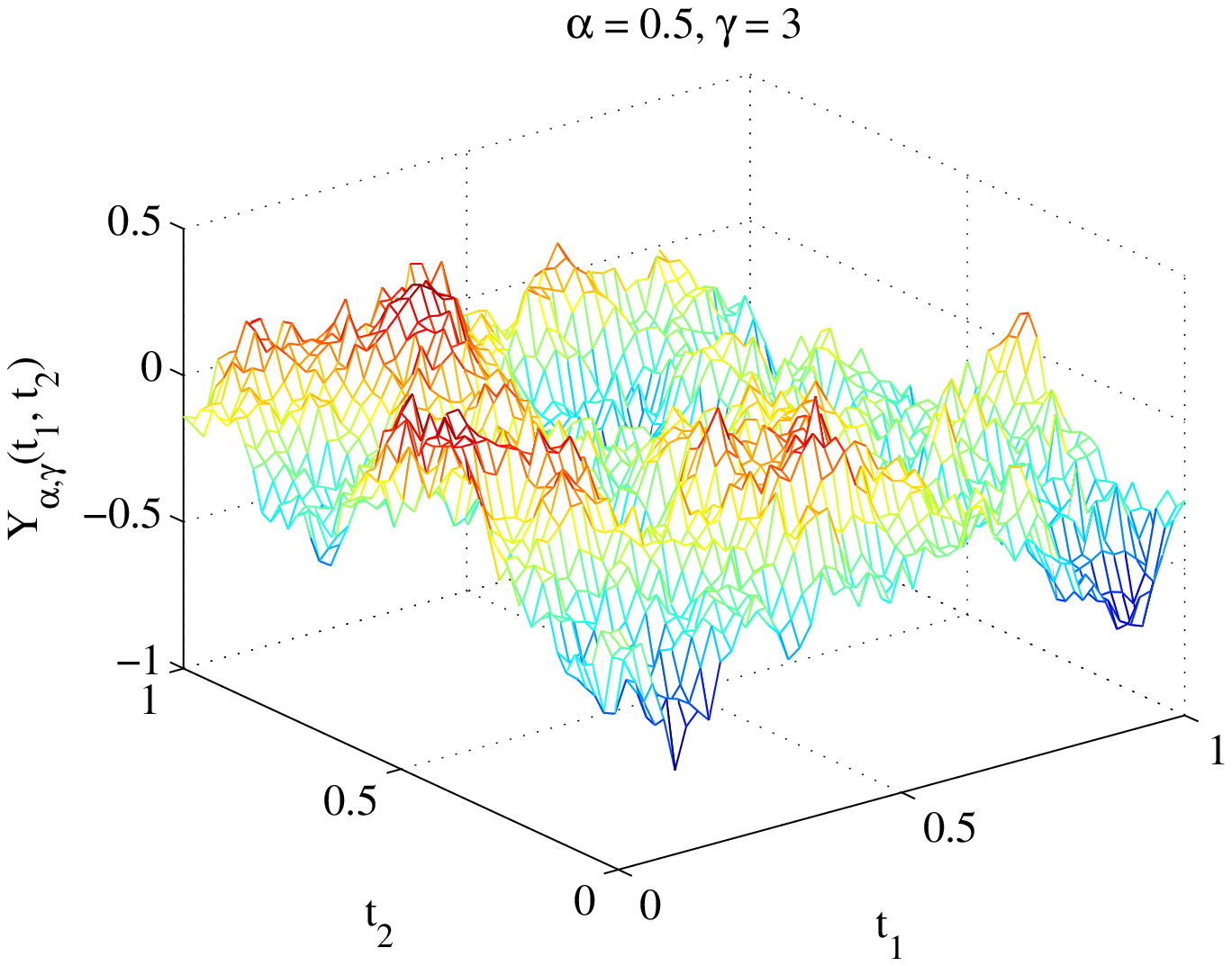}\caption{\label{Fig5} Simulations of the Whittle-Mat$\acute{\text{e}}$rn field $Y_{\alpha,\gamma}(t_1, t_2)$ for different values of $\alpha$ and $\gamma$. }\end{figure}
\begin{equation}\label{eq12} \begin{split}
C_{\alpha,\gamma}(\boldsymbol{t}-\boldsymbol{t}')=C_{\alpha,\gamma}(\boldsymbol{t},\boldsymbol{t}')=\langle
Y_{\alpha,\gamma}(\boldsymbol{t})Y_{\alpha,\gamma}(\boldsymbol{t}')\rangle =\frac{1}{(2\pi)^n}
\int\limits_{\mathbb{R}^n} \frac{e^{i\boldsymbol{\omega}.
(\boldsymbol{t}-\boldsymbol{t}')}}{\left(|\boldsymbol{\omega}|^{2\alpha}+\lambda^2\right)^{\gamma}}d^n\boldsymbol{\omega}.\end{split}
 \end{equation}
From this, we see that  $Y_{\alpha,\gamma}(\boldsymbol{t})$ is an isotropic field with spectral density
\begin{equation}\label{eq13}
S_{\alpha,\gamma}(\boldsymbol{\omega})=\frac{1}{(2\pi)^n}\frac{1}{\left(|\boldsymbol{\omega}|^{2\alpha}+\lambda^2\right)^{\gamma}}.
\end{equation}	 			
Note that in the case of $\alpha=1$, $\gamma>0$,  the field $Y_{1,\gamma}(\boldsymbol{t})$  is sometimes called Bessel field by some authors \cite{25, 26, 27, 28} based on the fact that the operator $\left(-\Delta +\lambda^2 \right)^{\gamma/2}$ is closely related to the Bessel potential with the WM covariance \eqref{eq1} equals to the Bessel kernel up to a multiplication constant \cite{22}. Since the function $\left(|\boldsymbol{\omega}|^{2\alpha}+\lambda^2\right)^{-\gamma/2}$  is in $L_2(\mathbb{R})$   if and only if $\alpha\gamma>n/2$, the field $Y_{\alpha,\gamma}(\boldsymbol{t})$  is only well-defined by \eqref{eq11} as an ordinary random field when $\alpha\gamma>n/2$. When $\alpha\gamma\leq n/2$, $Y_{\alpha,\gamma}(\boldsymbol{t})$  can be regarded as a generalized random field over the Schwarz space of test functions \cite{55}. In the following, when we study the properties of the random field $Y_{\alpha,\gamma}(\boldsymbol{t})$, we restrict to the case $\alpha\gamma>n/2$.
The two-dimensional GWM field with selected values of  $\alpha$ and $\gamma$   are simulated in Figure \ref{Fig5}. In next section we shall study the asymptotic properties of the covariance and the sample path properties of $Y_{\alpha,\gamma}(\boldsymbol{t})$.

Here we would like to remark that when $n=1$ and $\alpha=1$, the GWM process $Y_{1,\gamma}(t)$ is also called the Weyl fractional Ornstein--Uhlenbeck process or the Weyl fractional oscillator process \cite{nr8, nr9}, which can be considered as generalization of ordinary oscillator process driven by white noise.

\section{Asymptotic Properties of the Covariance Function}\label{sec3}

When $\alpha=1$, the covariance of GWM field \eqref{eq12} $C_{1,\gamma}(\boldsymbol{t})$ reduces to the WM class given by \eqref{eq1}.
However, \eqref{eq12} in general does not have a closed analytic form. It is interesting to note that the spectral density of the GWM field has the same functional form as both the characteristic function of generalized multivariate Linnik distribution \cite{29, 30} and the covariance function of generalized Cauchy class in  $\mathbb{R}^d$ \cite{31, 32}. Thus the covariance of the GWM field, the generalized multivariate Linnik distribution, and the spectral density of the random field belonging to the generalized Cauchy class all   should have the same analytic and asymptotic properties. These properties have been considered for the generalized Linnik distribution in $\mathbb{R}$, and the multivariate Linnik distribution for the special case with $\alpha\in (0,1)$  and $\gamma=1$, and for the spectral density of the random field of generalized Cauchy class by Kotz et al. \cite{29} and Ostrovskii \cite{30}, and Lim and Teo \cite{32} respectively. Thus the results obtained in \cite{29, 30, 32} can be translated directly to the covariance of the GWM field.

For general $\alpha$  and  $\gamma$, we can use a theorem of Bochner \cite{33} which says that being an isotropic covariance function, $C_{\alpha,\gamma}(\boldsymbol{t})$ has a spectral representation given by
 \begin{equation}\label{eq15}\begin{split}
C_{\alpha,\gamma}(\boldsymbol{t})=(2\pi)^{\frac{n}{2}}\int_0^{\infty}
\frac{J_{\frac{n-2}{2}}(\omega|\boldsymbol{t}|)}{(\omega|\boldsymbol{t}|)^{\frac{n-2}{2}}}
S_{\alpha,\gamma}(\omega)\omega^{n-1}d\omega\\=\frac{|\boldsymbol{t}|^{\frac{2-n}{2}}}{(2\pi)^{\frac{n}{2}}}\int_0^{\infty}
\frac{J_{\frac{n-2}{2}}( \omega|\boldsymbol{t}|)}{(\omega^{2\alpha}+\lambda^2)^{\gamma}}\omega^{\frac{n}{2}}d\omega.
\end{split}\end{equation}
Here $J_{\nu}(z)$  is the Bessel function of the first kind of order $\nu$. Now the result on a representation of the spectral density of the random field of generalized Cauchy class in \cite{32}  can be applied and we find that for $\alpha\in (0,1)$, the covariance function $C_{\alpha,\gamma}(\boldsymbol{t})$  has another representation given by
\begin{equation}\label{eq16}
C_{\alpha,\gamma}(\boldsymbol{t})=-\frac{|\boldsymbol{t}|^{\frac{2-n}{2}}}{2^{\frac{n-2}{2}}\pi^{\frac{n+2}{2}}}\text{Im}\,\int_0^{\infty}
\frac{K_{\frac{n-2}{2}}(u|\boldsymbol{t}|)}{\left(e^{i\pi\alpha}u^{2\alpha}+\lambda^2\right)^{\gamma}}u^{\frac{n}{2}}du.
\end{equation}
In fact, for all $\alpha\in (0,1)$  and $\gamma>0$, the integral in \eqref{eq16} is convergent when $\boldsymbol{t}\neq 0$. Together with \eqref{eq1}, we find that for $\alpha\gamma\leq n/2$, $Y_{\alpha,\gamma}(\boldsymbol{t})$    can be considered as a random field with infinite variance and with covariance given by \eqref{eq16} if $\alpha\in (0,1)$; and by \eqref{eq1} if $\alpha=1$. Since
\begin{equation*}
K_{\nu}(z)\sim \sqrt{\frac{\pi}{2z}}e^{-z}\hspace{1cm}\text{as}\;\;
z\rightarrow \infty,
\end{equation*}
we can use \eqref{eq16} to effectively calculate the numerical values of $C_{\alpha,\gamma}(\boldsymbol{t})$. On the other hand, we can also use \eqref{eq16} to study the large $|\boldsymbol{t}|$  behavior of the covariance function  $C_{\alpha,\gamma}(\boldsymbol{t})$ when $\alpha\in (0,1)$. More precisely, using the formula
\begin{equation*}
\frac{1}{(1+z)^{\gamma}}=\sum_{j=0}^{\infty}\frac{\Gamma(\gamma+j)}{\Gamma(\gamma)}\frac{(-1)^j}{j!}
z^j,
\end{equation*}
and the formula
\begin{equation*}
\int_0^{\infty} x^{\mu}K_{\nu}(x) dx
=2^{\mu-1}\Gamma\left(\frac{1+\mu+\nu}{2}\right)\Gamma\left(\frac{1+\mu-\nu}{2}\right),
\end{equation*}
(\cite{34}, \#6.561, no.16), we find that if $\alpha\in (0,1)$, then when $|\boldsymbol{t}|\rightarrow \infty$, we have
 \begin{equation}\label{eq17}\begin{split}
&C_{\alpha,\gamma}(\boldsymbol{t})=-\frac{|\boldsymbol{t}|^{-n}}{2^{\frac{n-2}{2}}\pi^{\frac{n+2}{2}}}\text{Im}\,\int_0^{\infty}
\frac{K_{\frac{n-2}{2}}(u)}{\left(e^{i\pi\alpha}\frac{u^{2\alpha}}{|\boldsymbol{t}|^{2\alpha}}
+\lambda^2\right)^{\gamma}}u^{\frac{n}{2}}du\\
\sim
&-\frac{|\boldsymbol{t}|^{-n}}{2^{\frac{n-2}{2}}\pi^{\frac{n+2}{2}}}\text{Im}\Biggl\{
\sum_{j=0}^{\infty}\frac{\Gamma(\gamma+j)}{\Gamma(\gamma)}\frac{(-1)^j}{j!}\lambda^{-2\gamma-j}
e^{i\pi\alpha j}|\boldsymbol{t}|^{-2\alpha j}\int_0^{\infty} u^{2\alpha j+
\frac{n}{2}}K_{\frac{n-2}{2}}(u)du \Biggr\}\\
\sim &\frac{1}{\pi^{\frac{n+2}{2}}}
\sum_{j=1}^{\infty}\frac{\Gamma(\gamma+j)}{\Gamma(\gamma)}\frac{(-1)^{j-1}}{j!}
\Gamma\left(\alpha j+1\right) \Gamma\left(\alpha j
+\frac{n}{2}\right)2^{2\alpha j}\lambda^{-2\gamma-j} \sin(\pi\alpha
j)|\boldsymbol{t}|^{-2\alpha j-n}.\end{split}
\end{equation}
In particular, the leading term of  $C_{\alpha,\gamma}(\boldsymbol{t})$ when  $|\boldsymbol{t}|\rightarrow \infty$  is
\begin{equation}\label{eq18}
C_{\alpha,\gamma}(\boldsymbol{t}) \sim
\frac{2^{2\alpha}\lambda^{-2\gamma-1}\gamma}{\pi^{\frac{n+2}{2}}}\Gamma(\alpha+1)\Gamma
\left(\alpha+\frac{n}{2}\right)\sin (\pi\alpha)|\boldsymbol{t}|^{-2\alpha-n}.
\end{equation}	 		
Note that the order of the leading term $|\boldsymbol{t}|^{-2\alpha-n}$  only depends on  $\alpha$. In other words, the large time asymptotic behavior of the covariance function varies as  $|\boldsymbol{t}|^{-2\alpha-n}$ and does not depend on $\gamma$. When $\alpha=1$, we cannot use \eqref{eq16}. However, we can obtain the large-$|\boldsymbol{t}|$  behavior of  $C_{1,\gamma}(\boldsymbol{t})$  from the explicit formula \eqref{eq1} and the asymptotic formula for $K_{\nu}(z)$   (\cite{34}, \#8.451, no.6) which give
\begin{equation}\label{eq19}\begin{split}
C_{1,\gamma}(\boldsymbol{t})\sim
\frac{2^{\frac{1-n}{2}-\gamma}}{\pi^{\frac{n-1}{2}}\Gamma(\gamma)}
e^{-\lambda|\boldsymbol{t}|}\sum_{j=0}^{\infty}\Biggl\{\frac{\Gamma\left(\gamma+j
-\frac{n-1}{2}\right)}{\Gamma\left(\gamma-j-\frac{n-1}{2}\right)} \frac{1}{2^j
j!} \lambda^{-j-\gamma+\frac{n-1}{2}} |\boldsymbol{t}|^{-j+\gamma-\frac{n+1}{2}}\Biggr\}.\end{split}
\end{equation}
Notice that in this case, $C_{1,\gamma}(\boldsymbol{t})$  decays exponentially and the leading term is
\begin{equation}\label{eq20}
C_{1,\gamma}(\boldsymbol{t})\sim\frac{2^{\frac{1-n}{2}-\gamma}\lambda^{-\gamma+\frac{n-1}{2}}}
{\pi^{\frac{n-1}{2}}\Gamma(\gamma)}
e^{-\lambda|\boldsymbol{t}|}|\boldsymbol{t}|^{\gamma-\frac{n+1}{2}}.
\end{equation}
\begin{figure}
\epsfxsize=0.49\linewidth \epsffile{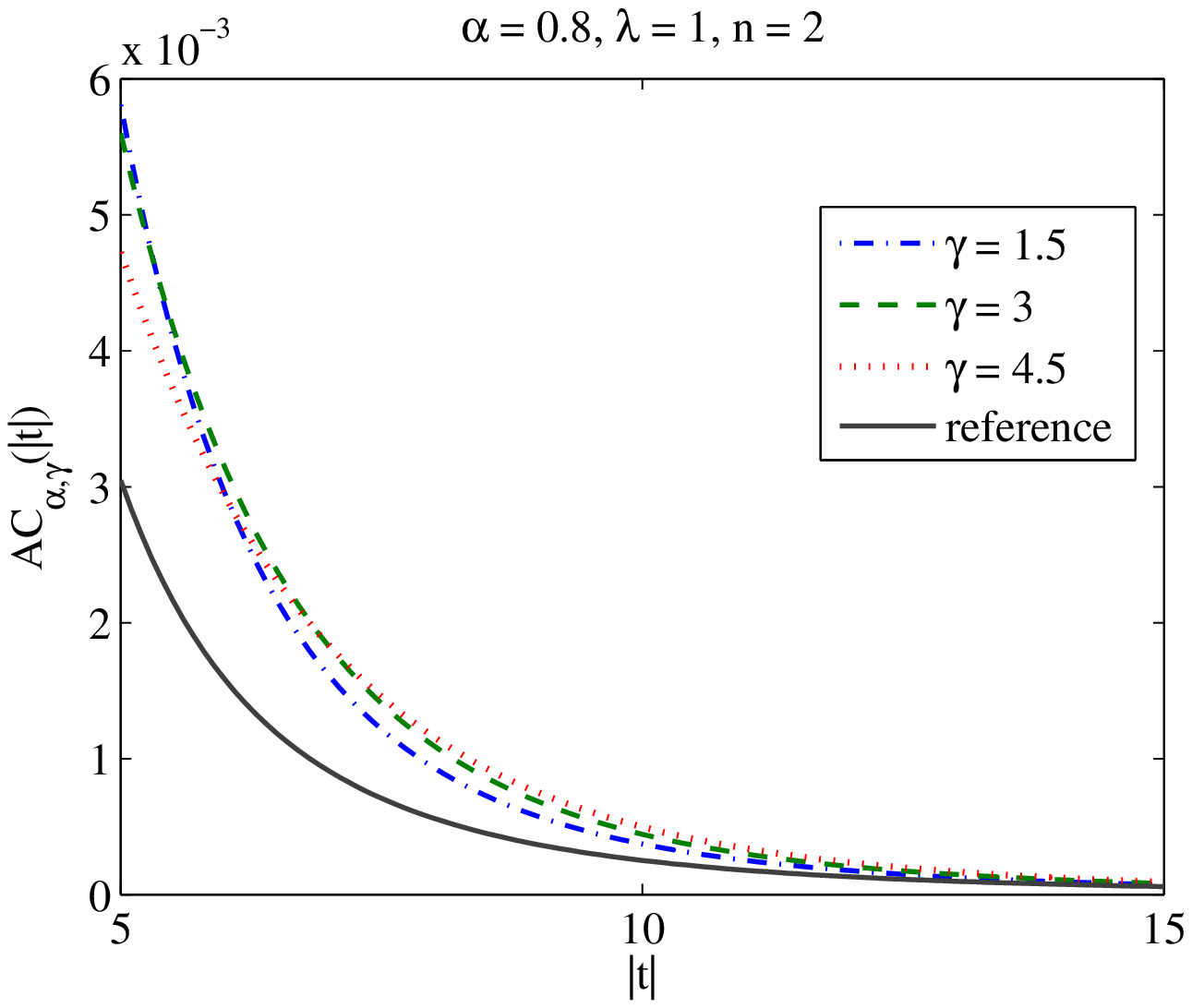}  \epsfxsize=0.49\linewidth
\epsffile{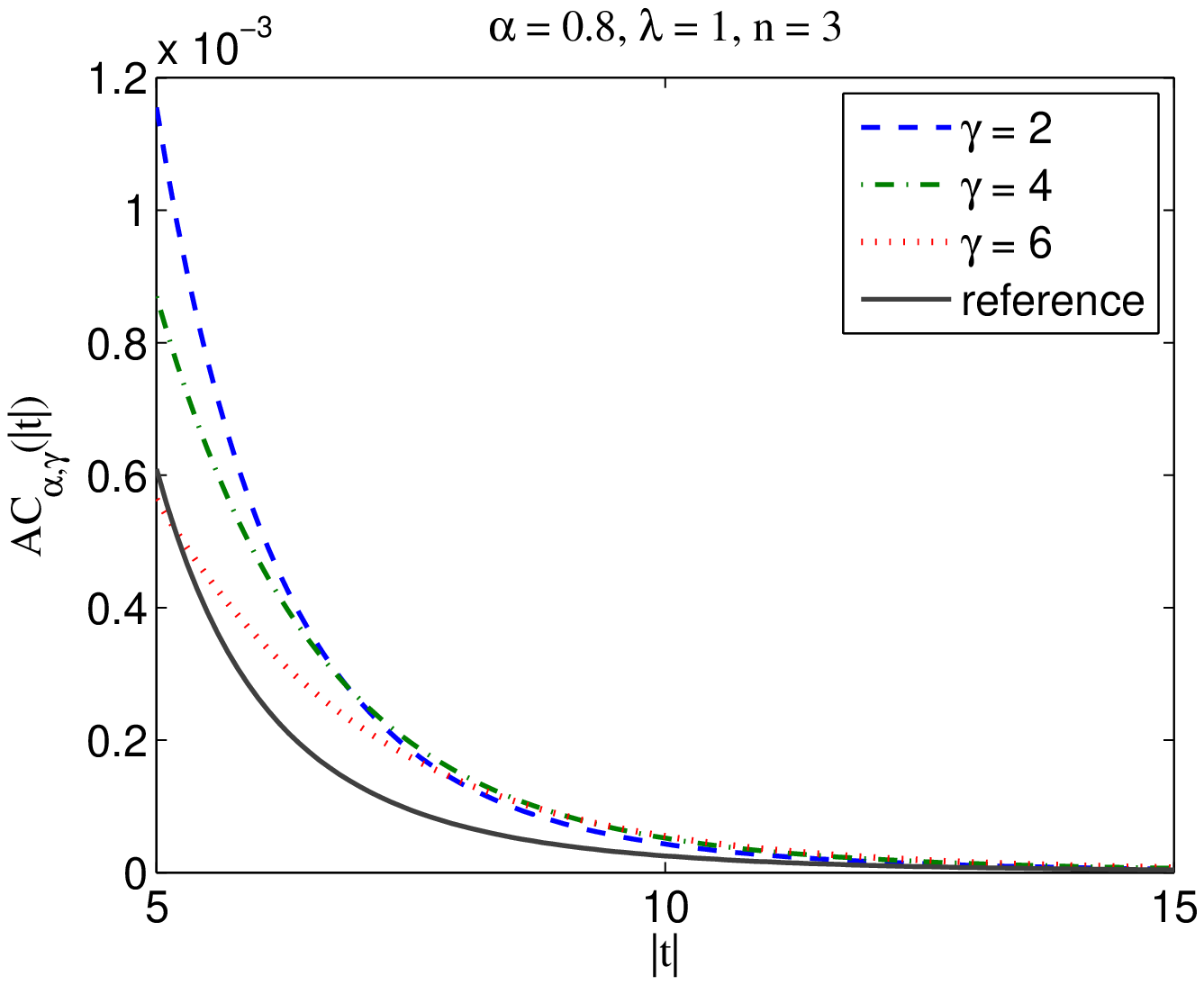}
 \caption{\label{Fig2}The graph of $AC_{\alpha,\gamma}(|\boldsymbol{t}|)$, $A=2^{2\alpha}\gamma \Gamma(\alpha+1)\Gamma(\alpha+n/2)\sin(\pi\alpha)/(\pi^{(n+2)/2}\lambda^{2\gamma+1})$, as   a function of $|\boldsymbol{t}|$. The reference curve is $y= |\boldsymbol{t}|^{-2\alpha -n}$. }
 \end{figure}
	
To study the local properties of the GWM field $Y_{\alpha,\gamma}(\boldsymbol{t})$ such as H\"older continuity, local asymptotic self similarity and Hausdorff dimension of the graph, we need to study the small-$|\boldsymbol{t}|$ behavior of the variogram
\begin{equation}\label{eq21}
\sigma_{\alpha,\gamma}^2(\boldsymbol{t}):=\langle
\left[Y_{\alpha,\gamma}(\boldsymbol{t})-Y_{\alpha,\gamma}(\mathbf{0})\right]^2\rangle
\end{equation}
of the increment field $Y_{\alpha,\gamma}(\boldsymbol{t})-Y_{\alpha,\gamma}(\mathbf{0})$. Notice that
\begin{equation}\label{eq22}
\sigma_{\alpha,\gamma}^2(\boldsymbol{t}) =
2(C_{\alpha,\gamma}(\mathbf{0})-C_{\alpha,\gamma}(\boldsymbol{t})),
\end{equation}
and the variance of  $Y_{\alpha,\gamma}(\boldsymbol{t})$  is given explicitly by
\begin{equation}\label{eq23}\begin{split}
\left\langle \left[Y_{\alpha,\gamma}(\boldsymbol{t})\right]^2\right\rangle
=&C_{\alpha,\gamma}(\mathbf{0})=\frac{1}{(2\pi)^n}\int\limits_{\mathbb{R}^n}\frac{1}{\left(|\boldsymbol{\omega}|^{2\alpha}+\lambda^2
\right)^{\gamma}}d^n\boldsymbol{\omega}\\
=&\frac{1}{2^{n-1}\pi^{\frac{n}{2}}\Gamma\left(\frac{n}{2}\right)}
\int_0^{\infty}\frac{\omega^{n-1}d\omega}{\left(\omega^{2\alpha}+\lambda^2\right)^{\gamma}}\\
=&\frac{\lambda^{\frac{n}{\alpha}-2\gamma}}{2^n
\pi^{\frac{n}{2}}\alpha
\Gamma\left(\frac{n}{2}\right)}\frac{\Gamma\left(\gamma-\frac{n}{2\alpha}\right)
\Gamma\left(\frac{n}{2\alpha}\right)}{\Gamma(\gamma)}.\end{split}\end{equation}
The $|\boldsymbol{t}|\rightarrow 0$  asymptotic properties of  $\sigma_{\alpha,\gamma}^2(\boldsymbol{t})$ depend on the arithmetic nature of  $\alpha$ and $\gamma$. To explore the leading behavior of  $\sigma_{\alpha,\gamma}^2(\boldsymbol{t})$  as  $|\boldsymbol{t}|\rightarrow 0$, we have to discuss the cases $\alpha\gamma\in\left( \frac{n}{2}, \frac{n+2}{2}\right)$, $\alpha\gamma=\frac{n+2}{2}$   and $\alpha\gamma>\frac{n+2}{2}$    separately. From \eqref{eq22} and \eqref{eq15},
 \begin{equation}\label{eq24}\begin{split}
\sigma^2_{\alpha,\gamma}(\boldsymbol{t})=\frac{-2}{(2\pi)^{\frac{n}{2}}}\int_0^{\infty}\left(\frac{J_{\frac{n-2}{2}}(k|\boldsymbol{t}|)}
{(k|\boldsymbol{t}|)^{\frac{n-2}{2}}}-\frac{1}{2^{\frac{n-2}{2}}\Gamma\left(\frac{n}{2}\right)}\right)
\frac{k^{n-1}}{(k^{2\alpha}+\lambda^2)^{\gamma}}dk.\end{split}
\end{equation}

\noindent
\textbf{Case I.} When  $\alpha\gamma\in\left( \frac{n}{2}, \frac{n+2}{2}\right)$, by making a change of variable  $k\mapsto k/|\boldsymbol{t}|$, \eqref{eq24} is transformed to
 \begin{equation}\label{eq25}\begin{split} \sigma_{\alpha,\gamma}^2(\boldsymbol{t})=\frac{-2|\boldsymbol{t}|^{2\alpha\gamma-n}}
{(2\pi)^{\frac{n}{2}}}\int_0^{\infty}\left(\frac{J_{\frac{n-2}{2}}(k)}
{k^{\frac{n-2}{2}}}-\frac{1}{2^{\frac{n-2}{2}}\Gamma\left(\frac{n}{2}\right)}\right)
\frac{k^{n-1}}{(k^{2\alpha}+\lambda^2|\boldsymbol{t}|^{2\alpha})^{\gamma}}dk.\end{split}\end{equation}
When  $|\boldsymbol{t}|\rightarrow 0$, the integral
 \begin{equation}\label{eq26}I(\boldsymbol{t})=\int_0^{\infty}\left(\frac{J_{\frac{n-2}{2}}(k)}
{k^{\frac{n-2}{2}}}-\frac{1}{2^{\frac{n-2}{2}}\Gamma\left(\frac{n}{2}\right)}\right)
\frac{k^{n-1}}{(k^{2\alpha}+\lambda^2|\boldsymbol{t}|^{2\alpha})^{\gamma}}dk\end{equation}
approaches a finite limit given by
\begin{equation*}
I=\int_0^{\infty}\left(\frac{J_{\frac{n-2}{2}}(k)}
{k^{\frac{n-2}{2}}}-\frac{1}{2^{\frac{n-2}{2}}\Gamma\left(\frac{n}{2}\right)}\right)
k^{n-2\alpha\gamma-1}dk.
\end{equation*}
Using regularization method (see appendix), it can be shown that
\begin{equation}\label{eq8_19_1}I=\frac{\Gamma\left(\frac{n}{2}-\alpha\gamma\right)}{2^{2\alpha\gamma-\frac{n}{2}}\Gamma(\alpha\gamma)}.\end{equation}
Therefore, as $|\boldsymbol{t}|\rightarrow 0$,
\begin{equation}\label{eq27}
\sigma_{\alpha,\gamma}^2(\boldsymbol{t})=-\frac{1}{2^{2\alpha\gamma-1}\pi^{\frac{n}{2}}}
\frac{\Gamma\left(\frac{n}{2}-\alpha\gamma\right)}{\Gamma(\alpha\gamma)}|\boldsymbol{t}|^{2\alpha\gamma-n}
+o(|\boldsymbol{t}|^{2\alpha\gamma-n}).
\end{equation}
Notice that when $\alpha\in \left(\frac{n}{2}, \frac{n+2}{2}\right)$, the leading order of $\sigma_{\alpha,\gamma}^2(\boldsymbol{t})$  depends on  $\alpha$ and  $\gamma$ only in the combination $\alpha\gamma$. By letting  $\gamma=\gamma'/\alpha$ gives  $\alpha\gamma=\gamma'$. Hence the   $|\boldsymbol{t}|\rightarrow 0$ asymptotic properties of  $\sigma_{\alpha,\gamma}^2(\boldsymbol{t})$ vary as  $|\boldsymbol{t}|^{2\gamma'-n}$ which is independent of $\alpha$.

\noindent
\textbf{Case II.} When  $\alpha\gamma>\frac{n+2}{2}$, using the fact that
\begin{equation*}
J_{\nu}(z)=\frac{z^{\nu}}{2^{\nu}}\sum_{j=0}^{\infty} \frac{(-1)^j
z^{2j}}{2^{2j}j!\Gamma(\nu+j+1)},
\end{equation*}	
(\cite{34},\#8.402), we find that as $|\boldsymbol{t}|\rightarrow 0$,\begin{equation}\label{eq28}
\frac{J_{\frac{n-2}{2}}(k|\boldsymbol{t}|)}
{(k|\boldsymbol{t}|)^{\frac{n-2}{2}}}-\frac{1}{2^{\frac{n-2}{2}}\Gamma\left(\frac{n}{2}\right)}=-
\frac{(k|\boldsymbol{t}|)^2}{2^{\frac{n+2}{2}}\Gamma\left(\frac{n+2}{2}\right)}+o(|\boldsymbol{t}|^2).
\end{equation}Therefore, \eqref{eq24} gives
\begin{equation}\label{eq29}\begin{split}
\sigma^2_{\alpha,\gamma}(\boldsymbol{t})=&\frac{|\boldsymbol{t}|^2}{2^n
\pi^{\frac{n}{2}}\Gamma\left(\frac{n+2}{2}\right)}\int_0^{\infty}
\frac{k^{n+1}dk}{\left(k^{2\alpha}+\lambda^2\right)^{\gamma}}+o(|\boldsymbol{t}|^2)\\
=&\frac{\lambda^{-2\gamma+\frac{n+2}{\alpha}}}{2^{n+1}
\pi^{\frac{n}{2}}\alpha\Gamma\left(\frac{n+2}{2}\right)}\frac{\Gamma\left(\gamma-\frac{n+2}{2\alpha}\right)
\Gamma\left(\frac{n+2}{2\alpha}\right)}{\Gamma(\gamma)}
|\boldsymbol{t}|^2+o(|\boldsymbol{t}|^2) \end{split}
\end{equation}as $|\boldsymbol{t}|\rightarrow 0$.

\begin{figure}
\epsfxsize=0.49\linewidth \epsffile{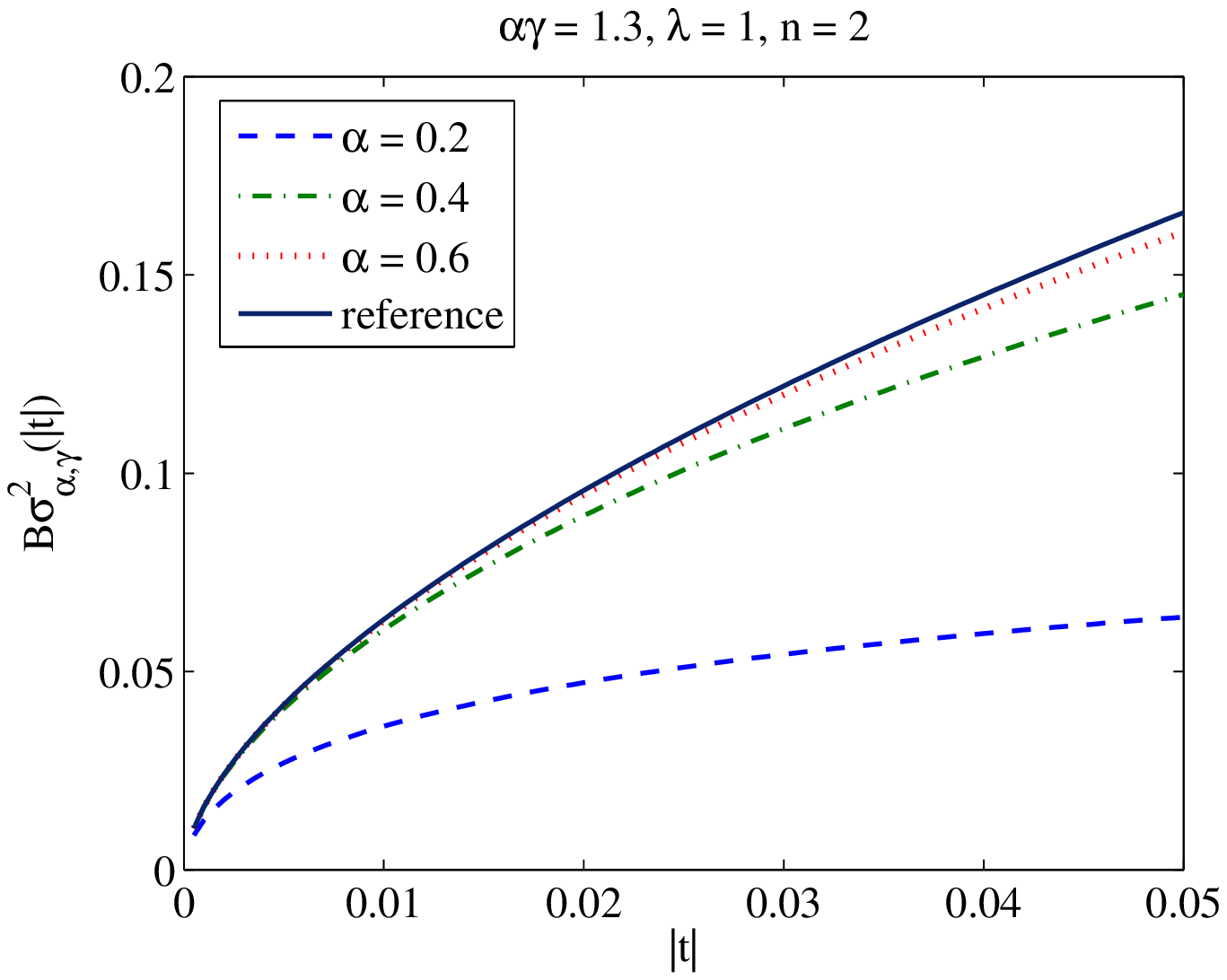}  \epsfxsize=0.49\linewidth
\epsffile{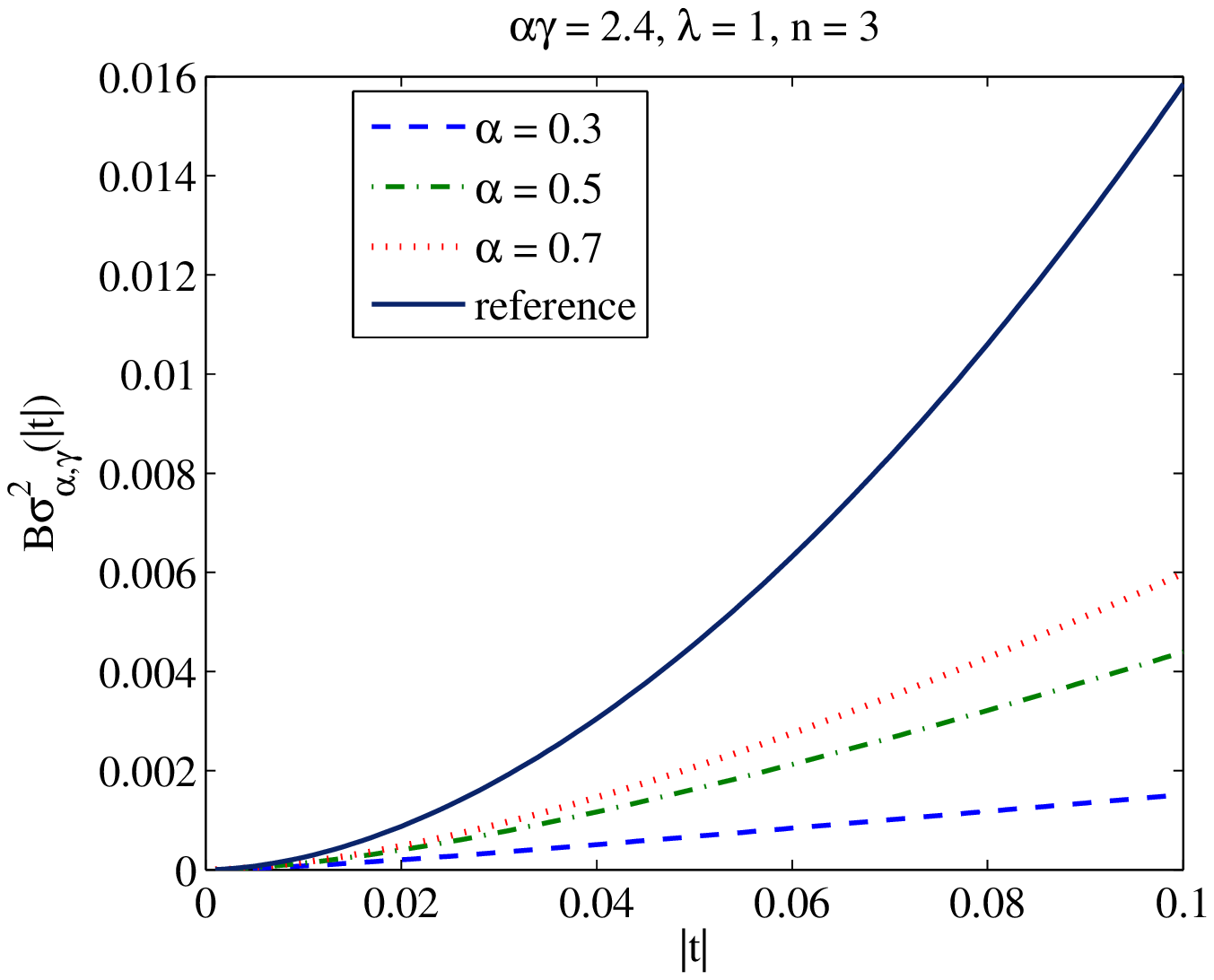}
 \caption{\label{Fig4}The graph of $B\sigma_{\alpha,\gamma}^2(\boldsymbol{t})$, $B=-\Gamma(n/2-\alpha\gamma)/(2^{2\alpha\gamma-1}\pi^{n/2}\Gamma(\alpha\gamma))$, as   a function of $|\boldsymbol{t}|$. The reference curve is $y=|\boldsymbol{t}|^{2\alpha\gamma-n}$.}
 \end{figure}	 			
\noindent
\textbf{Case III.} In the limiting case $\alpha\gamma=\frac{n+2}{2}$, eq. \eqref{eq25} gives
\begin{equation}\label{eq30}\begin{split}
\sigma_{\alpha,\gamma}^2(\boldsymbol{t})=\frac{-2|\boldsymbol{t}|^{2}}
{(2\pi)^{\frac{n}{2}}}\int_0^{\infty}\left(\frac{J_{\frac{n-2}{2}}(k)}
{k^{\frac{n-2}{2}}}-\frac{1}{2^{\frac{n-2}{2}}\Gamma\left(\frac{n}{2}\right)}\right)
\frac{k^{n-1}}{(k^{2\alpha}+\lambda^2|\boldsymbol{t}|^{2\alpha})^{\gamma}}dk.\end{split}
\end{equation}
However, now the integral
\begin{equation*}
I(\boldsymbol{t})=\int_0^{\infty}\left(\frac{J_{\frac{n-2}{2}}(k)}
{k^{\frac{n-2}{2}}}-\frac{1}{2^{\frac{n-2}{2}}\Gamma\left(\frac{n}{2}\right)}\right)
\frac{k^{n-1}}{(k^{2\alpha}+\lambda^2|\boldsymbol{t}|^{2\alpha})^{\gamma}}dk
\end{equation*}
does not have a finite limit when $ \boldsymbol{t}\rightarrow 0$. In the appendix, we show that
\begin{equation}
\label{eq31}I(\boldsymbol{t})=-\frac{1}{2^{\frac{n+2}{2}}\Gamma\left(\frac{n+2}{2}\right)}\log\frac{1}{|\boldsymbol{t}|}+A+o(1)\end{equation}
for some constant $A$. Therefore, as $|\boldsymbol{t}|\rightarrow 0$,
 \begin{equation}\label{eq32}
\sigma_{\alpha,\gamma}^2(\boldsymbol{t}) = \frac{1}{2^n
\pi^{\frac{n}{2}}\Gamma\left(\frac{n+2}{2}\right)}|\boldsymbol{t}|^2\log\frac{1}{|\boldsymbol{t}|}+O(|\boldsymbol{t}|^2).
\end{equation}
From \eqref{eq27}, \eqref{eq29} and \eqref{eq32}, we see that the behavior of the leading order term of  $\sigma_{\alpha,\gamma}^2(\boldsymbol{t})$  when $|\boldsymbol{t}|\rightarrow 0$  depends on $\gamma'-(n/2)$. If $\gamma'-(n/2)\in (0,1)$, the leading order term is of order $|\boldsymbol{t}|^{2\gamma'-n}$  which depends on the magnitude of $\gamma'-(n/2)$. If $\gamma'-(n/2)>1$, then the leading order term is of order $|\boldsymbol{t}|^2$, which loses dependence on $\gamma'-(n/2)$. In the borderline case  $\gamma'-(n/2)=1$, the leading order term is of order  $|\boldsymbol{t}|^2\log(1/|\boldsymbol{t}|)$.

	The graphs of $C_{\alpha,\gamma}(|\boldsymbol{t}|)$ when $|\boldsymbol{t}|$ is large and $\sigma_{\alpha,\gamma}^2(|\boldsymbol{t}|)$ when $|\boldsymbol{t}|$ is small   for some particular values of $\alpha$ and $\gamma$   are given in Figures  \ref{Fig2}  and \ref{Fig4} respectively.

 \section{Sample Path Properties of GWM Field   $\boldsymbol{Y_{\alpha,\gamma}}(\boldsymbol{t})$}\label{sec4}
Some basic sample path properties of the GWM field will be considered in this section.

\subsection{Continuity and differentiability}
When  $\alpha\gamma>n/2$ for which the field  $Y_{\alpha,\gamma}(\boldsymbol{t})$  is defined as an ordinary random field, the covariance function $C_{\alpha,\gamma}(\boldsymbol{t})$   \eqref{eq12} is continuous at  $\boldsymbol{t}=\mathbf{0}$. By a well-known result (see e.g.~\cite{35}), this implies that the field  $Y_{\alpha,\gamma}(\boldsymbol{t})$  is mean square (m.s.) continuous. One may then proceed to investigate the differentiability of the field  $Y_{\alpha,\gamma}(\boldsymbol{t})$. It turns out that $Y_{\alpha,\gamma}(\boldsymbol{t})$   is not always differentiable. In fact, a well-known result (see e.g.~\cite{35}) states that the m.s.~first partial derivative $\partial X(\boldsymbol{t})/\partial t_j$  of a stationary random field $X(\boldsymbol{t})$  exists if and only if the partial derivative  $\partial^2C(\boldsymbol{t})/\partial t_j^2$ exists at $\boldsymbol{t}=\mathbf{0}$, where $C(\boldsymbol{t})$  denotes the covariance function of $X(\boldsymbol{t})$. From our result in the previous section, we find that as $\boldsymbol{t} \rightarrow 0$,
\begin{equation}\label{eq33}\begin{split}
&C(\boldsymbol{t})-C(\boldsymbol{0})\\=&\begin{cases}B_1|\boldsymbol{t}|^{2\alpha\gamma-n}
+o(|\boldsymbol{t}|^{2\alpha\gamma-n}), \;\;&\text{if}\;\;
\alpha\gamma\in\left(\frac{n}{2},\frac{n+2}{2}\right)\\
B_2|\boldsymbol{t}|^2\log\frac{1}{|\boldsymbol{t}|}+O(|\boldsymbol{t}|^2), &\text{if}\;\;
\alpha\gamma=\frac{n+2}{2},\\
B_3|\boldsymbol{t}|^2+o(|\boldsymbol{t}|^2), &\text{if}\;\;\alpha\gamma>\frac{n+2}{2},
\end{cases}\end{split}
\end{equation}
for some constants $B_1, B_2, B_3$. It is easy to check that for the radial function $f(\boldsymbol{t})=|\boldsymbol{t}|^h$, the second partial derivative   $\partial^2f(\boldsymbol{t})/\partial t_j^2$ exists at $\boldsymbol{t}=\mathbf{0}$ if and only if $h\geq 2$. Therefore we conclude that the mean square partial derivatives of  $Y_{\alpha,\gamma}(\boldsymbol{t})$   exist if and only if $\alpha\gamma>\frac{n+2}{2}$, with a representation given by
\begin{equation}\label{eq34}
\frac{\partial Y_{\alpha,\gamma}}{\partial
t_j}(\boldsymbol{t})=\frac{i}{(2\pi)^{\frac{n}{2}}}\int\limits_{\mathbb{R}^n}\frac{\omega_j
e^{i\boldsymbol{t}.\boldsymbol{\omega}}\hat{\eta}(\boldsymbol{\omega})}{\left(|\boldsymbol{\omega}|^{2\alpha}+\lambda^2\right)^{\frac{\gamma}{2}}}
d^n\boldsymbol{\omega}.
\end{equation}
In fact, we can argue analogously that the m.s.~$j$-th order partial derivatives of   $Y_{\alpha,\gamma}(\boldsymbol{t})$  exist if and only if $\alpha\gamma>(n/2)+j$.

	For  $\alpha\gamma\in \left(\frac{n}{2}, \frac{n+2}{2}\right)$, the field   $Y_{\alpha,\gamma}(\boldsymbol{t})$  is not differentiable. Therefore, we would instead investigate the order of continuity of   $Y_{\alpha,\gamma}(\boldsymbol{t})$. Recall that a function $f$ is said to be H\"older continuous of order $h\in (0,1]$  if and only if
\begin{equation}\label{eq35}|f(\boldsymbol{t}')-f(\boldsymbol{t})|\leq K |\boldsymbol{t}'-\boldsymbol{t}|^h\hspace{1cm} \forall\;\; \boldsymbol{t}', \boldsymbol{t}\end{equation}  		
for some constant $K$. The sup of all $h$ where $f$ is H\"older continuous of order $h$ is called the H\"older exponent of $f$. For a centered isotropic  Gaussian random field $X(\boldsymbol{t})$, a concept of index-$\beta$  field was introduced by Adler \cite{35} which can be used to characterize the H\"older exponent of the sample paths of $X(\boldsymbol{t})$. More precisely, a theorem states that if $X(\boldsymbol{t})$ is an index-$\beta$  field, then with probability one, its sample paths have H\"older exponent equal to $\beta$, where $X(\boldsymbol{t})$ is called index-$\beta$  field if and only if \begin{equation}\label{eq36}\begin{split}
\beta=&\sup\left\{\tilde{\beta}\;:\;
\sigma(\boldsymbol{t})=o(|\boldsymbol{t}|^{\tilde{\beta}})\hspace{0.5cm}\text{as}\;\;
|\boldsymbol{t}|\rightarrow 0\right\}\\=&\inf\left\{\tilde{\beta}\;:\;
|\boldsymbol{t}|^{\tilde{\beta}}=o(\sigma(\boldsymbol{t}))\hspace{0.5cm}\text{as}\;\;
|\boldsymbol{t}|\rightarrow 0\right\}.\end{split}
\end{equation}
Here $\sigma(\boldsymbol{t})$ is defined as the square root of the variogram of $X(\boldsymbol{t})$, i.e., $\sigma(\boldsymbol{t})=\sqrt{\left\langle \left[X(\boldsymbol{t})-X(\boldsymbol{0})\right]^2\right\rangle}$. For the field  $Y_{\alpha,\gamma}(\boldsymbol{t})$  we are considering,  it is immediate to conclude that from \eqref{eq27}, \eqref{eq29} and \eqref{eq32} that if $\alpha\gamma\in \left(\frac{n}{2}, \frac{n+2}{2}\right)$, then  $Y_{\alpha,\gamma}(\boldsymbol{t})$  is an indexed  $\left(\alpha\gamma-(n/2)\right)$ field; whereas if $\alpha\gamma\geq \frac{n+2}{2}$, then $Y_{\alpha,\gamma}(\boldsymbol{t})$   is an index-1 field. Therefore, we have for  $\alpha\gamma\in \left(\frac{n}{2}, \frac{n+2}{2}\right)$, the sample paths of $Y_{\alpha,\gamma}(\boldsymbol{t})$  is H\"older  continuous of order $\alpha\gamma-(n/2)$  with probability one. For $\alpha\gamma> \frac{n+2}{2}$, it can be shown by considering the gradient field $\nabla Y_{\alpha,\gamma}(\boldsymbol{t})=\left( \partial Y_{\alpha,\gamma}(\boldsymbol{t})/\partial t_1, \ldots,  \partial Y_{\alpha,\gamma}(\boldsymbol{t})/\partial t_n\right)$  that the sample paths of  $Y_{\alpha,\gamma}(\boldsymbol{t})$ are differentiable.

\subsection{Fractal dimension}
For a non-differentiable function $f$, ordinary definition of dimension, which is always a nonnegative integer, is inadequate to measure the dimensionality of its image or graph. A more appropriate definition of dimension is called fractal of Hausdorff dimension which can be any non-negative real number. Definition and basic properties of fractal dimension can be obtained in the book \cite{36}. Here we would like to make use of the following result. For an index-$\beta$  field in  $\mathbb{R}^n$, with probability one, the fractal dimension of the image and graph of its sample path are 1 and  $n+1-\beta$  respectively. Thus, with probability one, the image of the sample path of $Y_{\alpha,\gamma}(\boldsymbol{t})$   always has fractal dimension one. The result is more interesting for the fractal dimension of the graphs. If  $\alpha\gamma\in \left(\frac{n}{2}, \frac{n+2}{2}\right)$, then with probability one, the graph of the sample path of  $Y_{\alpha,\gamma}(\boldsymbol{t})$ has dimension  $\frac{3n}{2}+1-\alpha\gamma$, a real number between $n$ and $n+1$. However, when $\alpha\gamma$  exceeds the point $\frac{n+2}{2}$, then with probability one, the graph of the sample path of  $Y_{\alpha,\gamma}(\boldsymbol{t})$  always has dimension equal to $n$. This is reasonable since the sample path   of $Y_{\alpha,\gamma}(\boldsymbol{t})$ becomes differentiable when $\alpha\gamma > \frac{n+2}{2}$. In fact, Figure \ref{Fig5} show clearly that the fractal dimension of the graph of $Y_{\alpha,\gamma}(t_1, t_2)$ depend on $\gamma'=\alpha\gamma$.

\subsection{ Local self-similarity}
Self-similarity is an important property of fractals. Intuitively, a field is called self-similar if it is invariant under appropriate scaling. For a random field $X(\boldsymbol{t})$, we say that it is self-similar of order $H$ if and only if for any $r>0$, the law of  the field $X(r\boldsymbol{t})$  is the same as the law of the field $r^{H}X(\boldsymbol{t})$. It is well-known that a stationary random field  cannot be self-similar \cite{37}.  In fact, up to a constant multiplicative factor, the only $H$-self-similar centered Gaussian random field with stationary increments is the fractional L$\acute{\text{e}}$vy Brownian field $B_H(\boldsymbol{t})$ of index $H$ with covariance
\begin{equation}\label{eq8_20_2}
\left\langle B_H(\boldsymbol{s})B_H(\boldsymbol{t})\right\rangle = \frac{1}{2}\left(
|\boldsymbol{t}|^{2H}+|\boldsymbol{s}|^{2H}-|\boldsymbol{t}-\boldsymbol{s}|^{2H}\right).
\end{equation}	 				
This excludes the possibility for   $Y_{\alpha,\gamma}(\boldsymbol{t})$  being a self-similar random field. However,   $Y_{\alpha,\gamma}(\boldsymbol{t})$ satisfies a weaker self-similar property known as local self-similarity considered by Kent and Wood \cite{38}. A centered stationary Gaussian field is locally self-similar of order  $\beta/2$ if its covariance  $C(\boldsymbol{t})$ satisfies for $| \boldsymbol{t}|\rightarrow 0$,\begin{equation}\label{eq38}
C(\boldsymbol{t})=C(\mathbf{0})-A|\boldsymbol{t}|^{\beta}\left[1+O(|\boldsymbol{t}|^{\delta})\right]
\end{equation}with $A>0$ and $\delta>0$.
The proof of eq. \eqref{eq27} shows that for $\alpha\gamma\in\left(\frac{n}{2},\frac{n+2}{2}\right)$,
\begin{equation*}
C_{\alpha,\gamma}(\boldsymbol{t}) =C_{\alpha,\gamma}(\mathbf{0})-A|\boldsymbol{t}|^{2\alpha\gamma-n}+o(|\boldsymbol{t}|^{2\alpha\gamma-n+\delta})
\end{equation*}
with \begin{equation}\label{eq8_20_1}A=-\frac{1}{2^{2\alpha\gamma}\pi^{\frac{n}{2}}}\frac{\Gamma\left(\frac{n}{2}-\alpha\gamma\right)}
{\Gamma(\alpha\gamma)}.\end{equation}  Hence $Y_{\alpha,\gamma}(\boldsymbol{t})$ is locally self-similar of order $\alpha\gamma-(n/2)$.

	There exists an equivalent way of characterizing self-similarity at a local scale called local asymptotical self-similarity  which was first introduced for multifractional Brownian motion \cite{39}. Recall that a random field $X(\boldsymbol{t})$ is called locally asymptotically self-similar with parameter $H\in (0,1)$  at a point $\boldsymbol{t}_0$ if the limit random field
\begin{equation}\label{eq40}
\left\{ T_{\boldsymbol{t}_0}(\boldsymbol{u})=\lim_{\rho\rightarrow 0^+}\frac{X(\boldsymbol{t}_0+\rho \boldsymbol{u})
-X(\boldsymbol{t}_0)}{\rho^{H}},\hspace{0.5cm}\boldsymbol{u}\in\mathbb{R}^n\right\}
\end{equation}exists and is nontrivial \cite{39}. In this case,  $T_{\boldsymbol{t}_0}(\boldsymbol{u})$ is called the tangent field of $X(\boldsymbol{t})$ at  $\boldsymbol{t}_0$. It can be directly verified that
\begin{equation}\label{eq41}\begin{split}
&\Bigl\langle \left[Y_{\alpha,\gamma}(\boldsymbol{t}_0+\rho \boldsymbol{u}) -Y_{\alpha,\gamma}(\boldsymbol{t}_0)\right]\left[Y_{\alpha,\gamma}(\boldsymbol{t}_0+\rho \boldsymbol{v})
-Y_{\alpha,\gamma}(\boldsymbol{t}_0)\right]\Bigr\rangle\\=&\frac{1}{2}\left(\sigma^2_{\alpha,\gamma}(\rho
\boldsymbol{u})+\sigma^2_{\alpha,\gamma}(\rho
\boldsymbol{v})-\sigma_{\alpha,\gamma}^2(\rho(\boldsymbol{u}-\boldsymbol{v}))\right).\end{split}
\end{equation}
Eq. \eqref{eq27} then shows that for  $\alpha\gamma\in\left(\frac{n}{2},\frac{n+2}{2}\right)$,
\begin{equation}\label{eq42}\begin{split}
&\lim_{\rho\rightarrow 0^+}\left\langle \frac{Y(\boldsymbol{t}_0+\rho \boldsymbol{u})
-Y(\boldsymbol{t}_0)}{\rho^{\alpha\gamma-\frac{n}{2}}}\frac{Y(\boldsymbol{t}_0+\rho \boldsymbol{v})
-Y(\boldsymbol{t}_0)}{\rho^{\alpha\gamma-\frac{n}{2}}}\right\rangle\\
=&A\left(|\boldsymbol{u}|^{2\alpha\gamma-n}
+|\boldsymbol{v}|^{2\alpha\gamma-n}-|\boldsymbol{u}-\boldsymbol{v}|^{2\alpha\gamma-n}\right).\end{split}
\end{equation}
Therefore, for  $\alpha\gamma\in\left(\frac{n}{2},\frac{n+2}{2}\right)$, $Y_{\alpha,\gamma}(\boldsymbol{t})$  is locally asymptotically self-similar with order $\alpha\gamma-(n/2)$. Its tangent field at a point $\boldsymbol{t}_0\in \mathbb{R}^n$  is independent of $\boldsymbol{t}_0$, and up to the multiplicative factor $2A$, the tangent is given by the fractional  L$\acute{\text{e}}$vy Brownian field  $B_H(\boldsymbol{u})$ \eqref{eq8_20_2} of order $\alpha\gamma-(n/2)$.

	From the results on H\"older exponent, fractal dimension and local self-similarity, it is found that they all depend on the parameters $\alpha$ and $\gamma$ in the combination $\gamma'=\alpha\gamma$. They are related to each other in such a way that if the H\"older exponent is $H=\min \{ \gamma'-(n/2), 1\}$, then the fractal dimension is $n+1-H$     and the order of local self-similarity is $H$ again if  $H <1$.

\subsection{Short Range Dependence}
Recall that a stationary random field $X(\boldsymbol{t})$ is said to have short range dependence (or short
memory) if the absolute value of its covariance function  $C(\boldsymbol{t})$ is integrable over  $\mathbb{R}^n$, that is
\begin{align*}
\int\limits_{\mathbb{R}_+^n} |C(\boldsymbol{t})|d^n\boldsymbol{t}<\infty.
\end{align*}By our result on the large $|\boldsymbol{t}|$ asymptotic behavior of
$C_{\alpha,\gamma}(\boldsymbol{t})$ \eqref{eq18} and \eqref{eq20}, we find
that when $|\boldsymbol{t}|\rightarrow \infty$,
$C_{\alpha,\gamma}(\boldsymbol{t})\sim |\boldsymbol{t}|^{-2\alpha-n}$ if $\alpha\in (0,1)$ and
$C_{\alpha,\gamma}(\boldsymbol{t})\sim e^{-\lambda|\boldsymbol{t}|}
|\boldsymbol{t}|^{\gamma-\frac{n+1}{2}} $ if $\alpha=1$. Using polar
coordinates, it can be verified easily that
$$\int\limits_{\boldsymbol{t}\in \mathbb{R}^n_+, |\boldsymbol{t}|>1}|\boldsymbol{t}|^{p}d^n\boldsymbol{t}<\infty$$ if and only if $p<-n$.
This immediately implies that $Y_{\alpha,\gamma}(\boldsymbol{t})$ has short range
dependence for all $\alpha$ and $\gamma$. Moreover, for $\alpha\in
(0,1)$, the short memory exponent $2\alpha+n$  depends only on
$\alpha$ and not on $\gamma$. Together with the result on local
properties such as H\"older exponent, fractal dimension and local
asymptotic self-similarity, this implies that the short range
dependence property and local properties of the GWM field $Y_{\alpha,\gamma}(\boldsymbol{t})$ are characterized separately
by  $\alpha$ and $\gamma'=\alpha\gamma$.  This should be compared to the random field with generalized Cauchy covariance \cite{31, 32}, which is an isotropic  random field with two parameters
that enables separate characterizations of long range dependence and
fractal dimension.

	Here we would like to remark on the Markov property for WM field and GWM field. In the case of Whittle field ($\alpha=\gamma=1$) in $n=2$ dimension, this problem has been studied by Pitt and Robeva \cite{40}, who showed that under certain technical conditions, the sharp Markov property is satisfied. They generalized the result to WM field (which they called Bessel field), and they verified that under some technical conditions the sharp Markov property holds for WM field with  $n+1/2 <\gamma < n+1$, $n\geq 1$ \cite{26, 27}. It will be interesting to see whether the arguments of Pitt and Robeva can be extended to the GWM field.

\section{ Application to wind speed modeling}WM field has been widely used in modeling \cite{10, 11, 12, 13, 41} geostastical data such as sea beam data, temperature, wind speed and soil data. In this section, we   show that the GWM process can be used to provide an alternative model for wind speed.

\begin{figure}
\epsfxsize=0.6\linewidth \epsffile{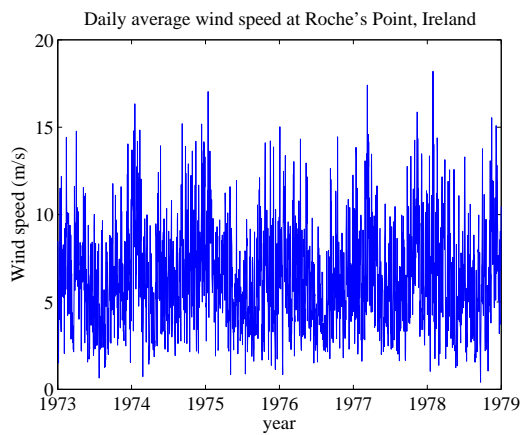} \caption{\label{Fig9}The daily average wind speed at Roche's Point, Ireland from 1973 to 1978.}\end{figure}

We analyze the average daily wind speed of Roche's Point in Ireland from 1973 to 1978 which consists of $N=365\times 6=2190$  data points\footnote{The data is obtained from Statlib (http://lib.stat.cmu.edu/datasets/) with the value for 29th, February, 1976  omitted.} (see Figure \ref{Fig9}).  The Irish wind data of 12 meteorological sites from 1961 to 1978 has been analyzed by several authors \cite{n2, n3, n4, n5} where they were more concerned with the  spatial correlation between the sites.  On the other hand, the von K$\acute{\text{a}}$rm$\acute{\text{a}}$n wind turbulence model \cite{16, 17} proposed a model with spectral density having the same functional form as the spectral density of the WM process.

\begin{figure}
\epsfxsize=0.6\linewidth \epsffile{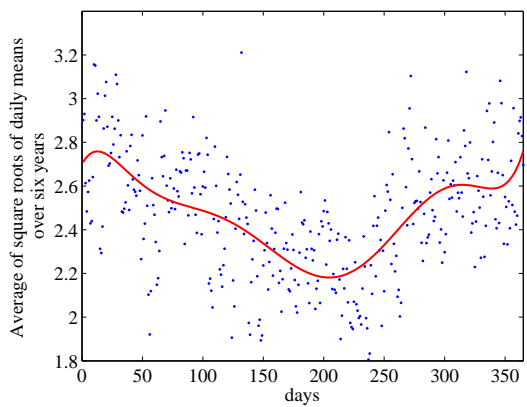} \caption{\label{Fig10}The average of the square roots of daily means over the 6 years for each day of the year and the fitted polynomial of degree 8.}\end{figure}

\begin{figure}
\epsfxsize=0.6\linewidth \epsffile{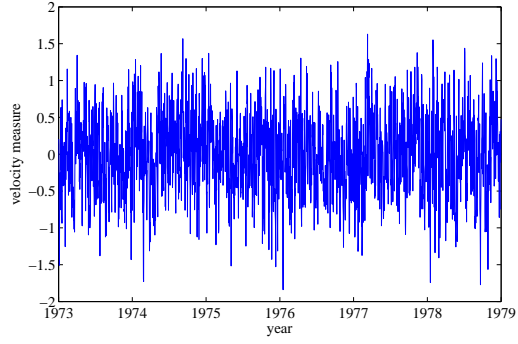} \caption{\label{Fig11}The deseasonalized data (velocity measures).}\end{figure}

\begin{figure}
\epsfxsize=0.49\linewidth \epsffile{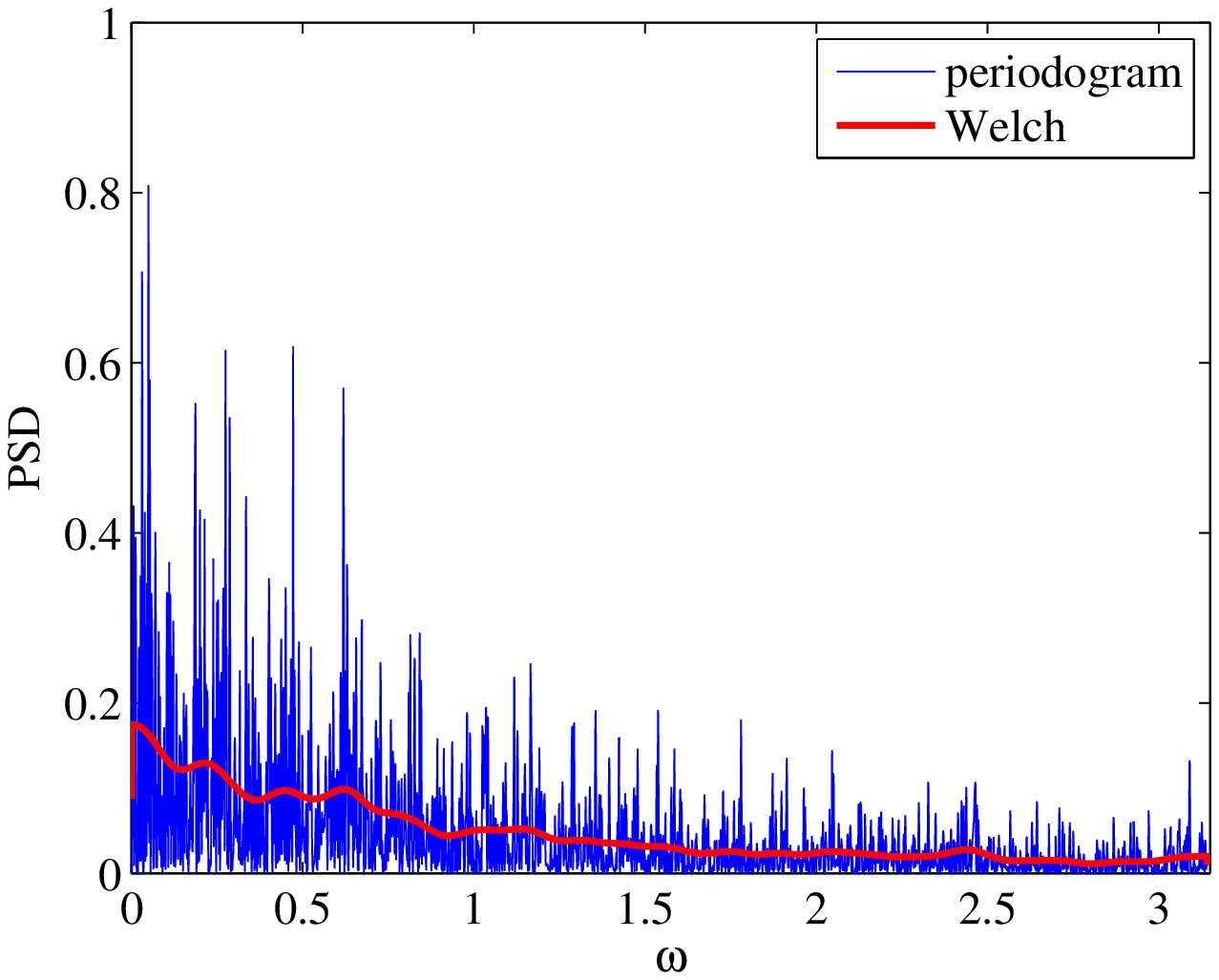} \epsfxsize=0.49\linewidth \epsffile{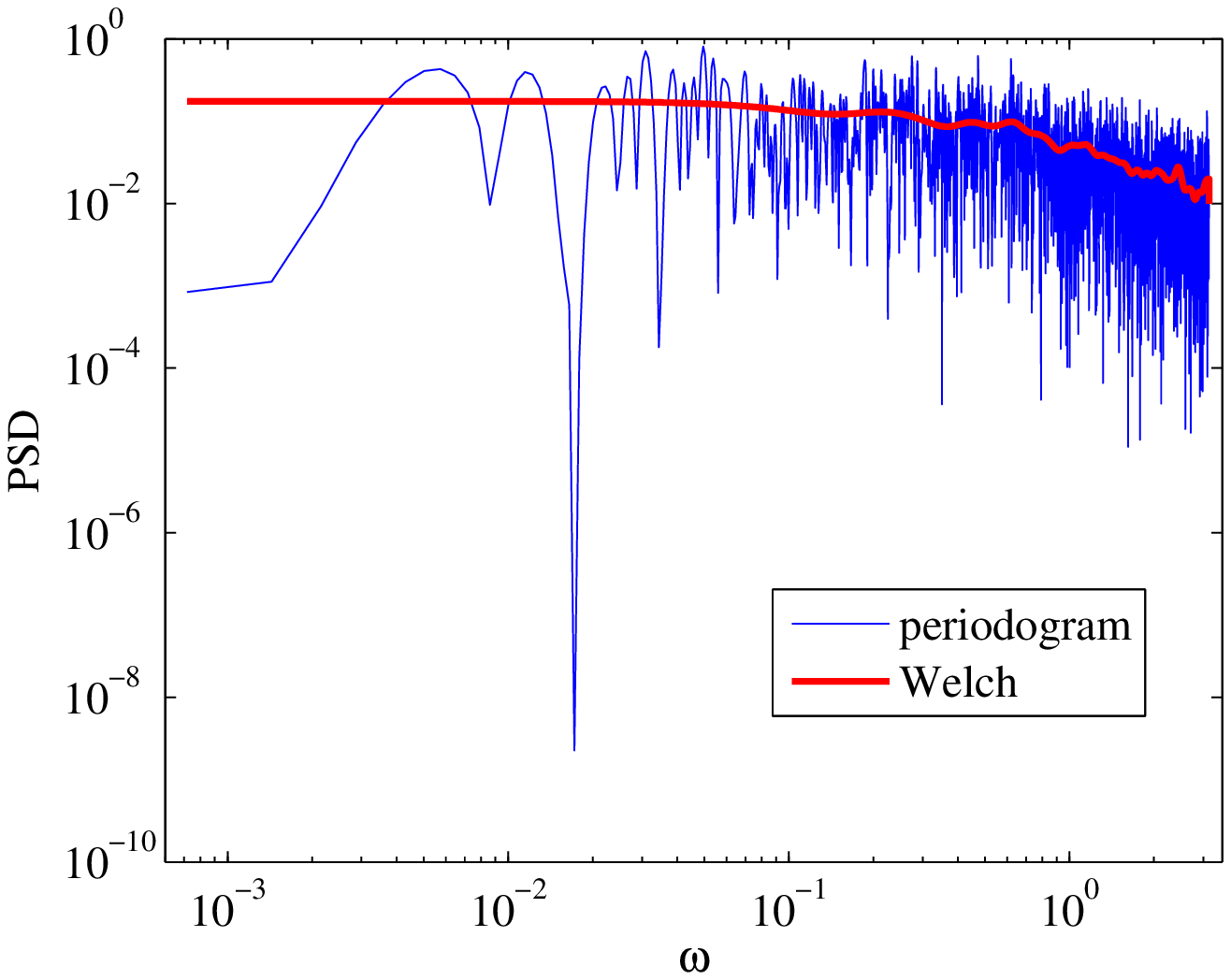} \caption{\label{Fig12}The estimated power spectral density of the velocity measures and the corresponding log-log-plot using (a) periodogram method, (b) Welch's method.}\end{figure}

As in \cite{n2}, we consider the seasonal effect by calculating the average of the square roots of the daily means over the 6 years for each day of the year, and then regressing the result with a polynomial of degree 8 (see Figure \ref{Fig10}). The deseasonalized data (Figure \ref{Fig11}) is obtained by subtracting the fitted polynomial from the square roots of daily means. It has zero mean and is referred to as the velocity measures. To justify that the velocity measures is a short memory process, we use the fact that a process has long memory if and only if it's spectral density diverges at $\omega=0$. For a discrete stationary random process $X_t$, $t=1,2,3,\ldots$, with covariance $C(t)$, an analog of spectral density is the power spectral density (PSD) defined by
\begin{align*}
\text{PSD}(\omega) = \frac{1}{2\pi}\sum_{j=-\infty}^{\infty} C(j)e^{-ij\omega}.
\end{align*}
It is a periodic function with period $2\pi$ and $S(2\pi -\omega)=S(\omega)$. If $X_t$ has an underlying continuous process $X(t)$ with spectral density $S(\omega)$ so that $X_t=X(t)$ when $t=1, 2, 3, \ldots$, then
\begin{align*}
\text{PSD}(\omega) =\sum_{j=-\infty}^{\infty}S(\omega-2\pi j).
\end{align*}There are  different  ways to estimate the power spectral density from a given sample $x_t, t=1,2,\ldots, N$ of $X_t$. One way is to use the periodogram method, where the estimate of $\text{PSD}(\omega)$ is given by the periodogram
\begin{equation*}
\frac{1}{2\pi N} \left|\sum_{j=1}^{N} x_j e^{-i\omega j}\right|^2.
\end{equation*}However, this method usually leads to   large fluctuations.  A better method which gives a smoother estimate is introduced by Welch \cite{n6} and improved by others (see e.g.~\cite{n10}). We estimate the PSD of the velocity measures using Welch's method by segmenting the data into 50\% overlapping blocks of length 73 and applying the Hamming window to each block. The resulting estimate for PSD is compared to the periodogram estimate in Figure \ref{Fig12}.
The behavior of PSD at $\omega\sim 0$ gives strong evidence that the velocity measures has short memory.

\begin{figure}
\epsfxsize=0.49\linewidth \epsffile{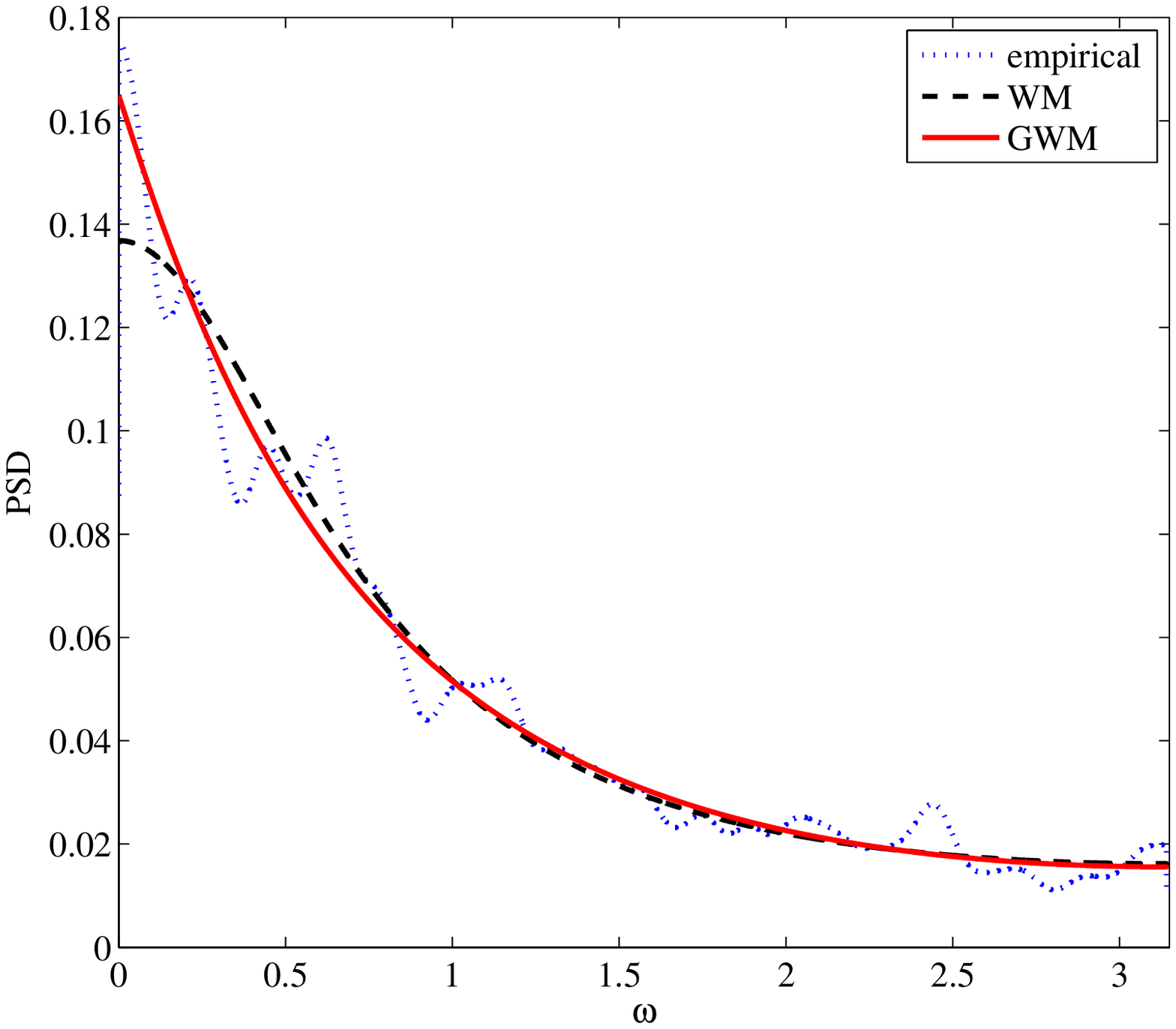} \epsfxsize=0.49\linewidth \epsffile{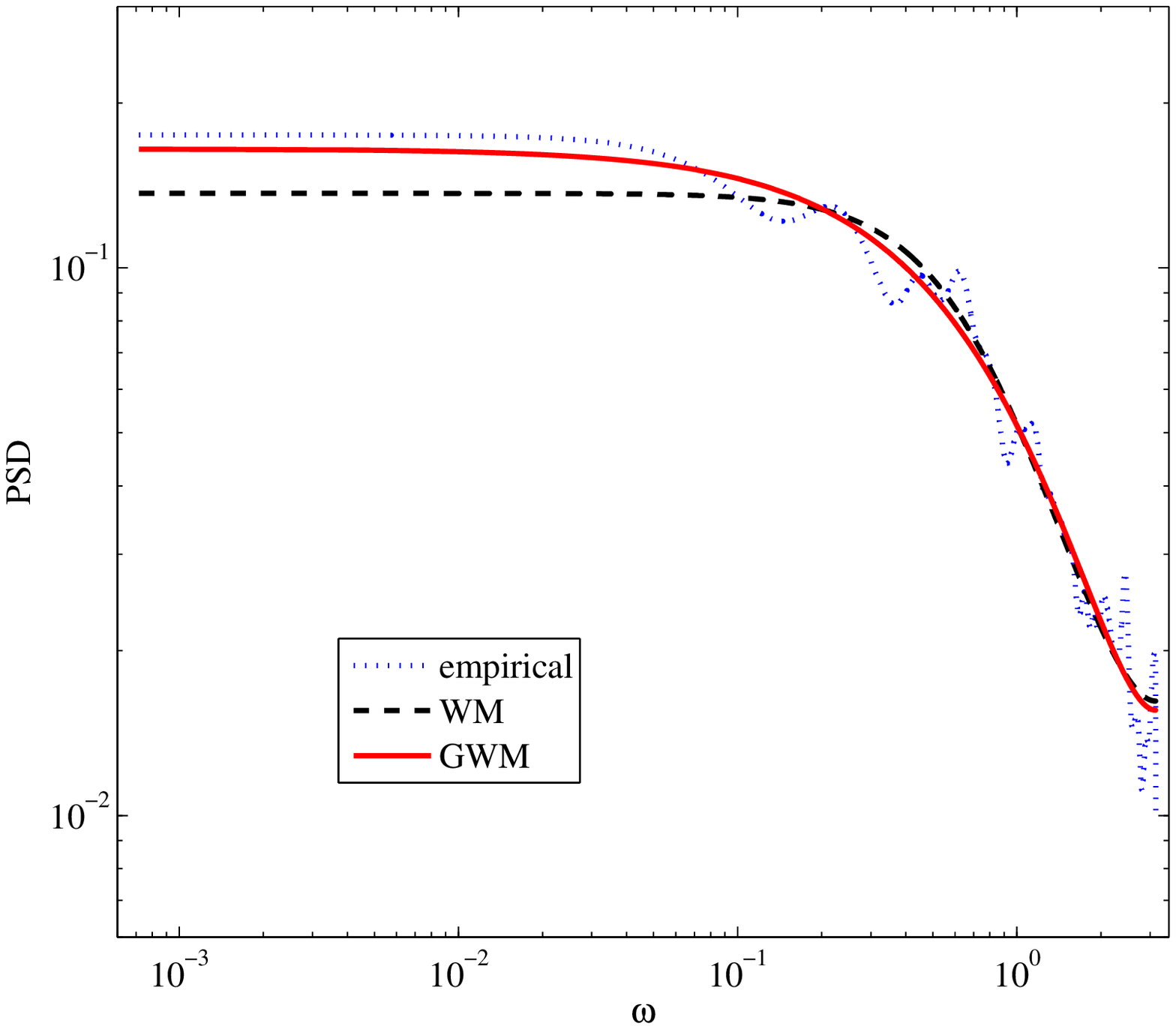} \caption{\label{Fig13}The empirical power spectral density of the velocity measures compared to the PSD of the WM model and the GWM model with parameters given in Table \ref{Tab1}, and the corresponding log-log-plot.}\end{figure}

In order to apply the GWM process to model the velocity measures, we consider the more general process parametrized by four parameters $\alpha, \gamma, K, \ell$:
\begin{equation*}
\hat{Y}_{\alpha,\gamma}^{K, \ell}(t) = K Y_{\alpha,\gamma}(\ell t)_{\lambda=1}=\frac{K}{\sqrt{2\pi}}\int_{-\infty}^{\infty}\frac{e^{i\ell t\omega}\hat{\eta}(\omega)}
{\left(|\omega|^{2\alpha}+1\right)^{\frac{\gamma}{2}}}d\omega,
\end{equation*}which has spectral density
\begin{equation*}
\hat{S}_{\alpha,\gamma}^{K, \ell}(\omega) = \frac{1}{2\pi} \frac{K\ell^{2\alpha\gamma-1}}{\left(|\omega|^{2\alpha}+\ell^{2\alpha}\right)^{\gamma}}=
K\ell^{2\alpha\gamma-1}S_{\alpha,\gamma}(\omega)_{\lambda=\ell^{2\alpha}}
\end{equation*}and covariance function
\begin{equation*}
\hat{C}_{\alpha,\gamma}^{K, \ell}(t)=K^2C_{\alpha,\gamma}(\ell t)_{\lambda=1}.
\end{equation*}
Notice that $\ell$ rescale the time parameter and $K$ rescale the data.
We need to determine the parameters $\alpha,\gamma, K, \ell$ so that the process $\hat{Y}_{\alpha,\gamma}^{K, \ell}(t)$ gives the best model to the velocity measures. For this purpose we use the maximum likelihood estimation (MLE) method strongly recommended by Stein (see e.g.~\cite{14}). Let $\boldsymbol{\Gamma}\!\left(\boldsymbol{\theta}\right)=\Gamma(\alpha,\gamma, K,\ell)$ be the covariance matrix $(\hat{C}_{\alpha,\gamma}^{K, \ell}(i-j))_{i, j=1}^{N}$. Since we assume that the velocity measure is a Gaussian process, the probability density function for $\boldsymbol{y}=(y_1, \ldots, y_N)^T$ having mean $\mathbf{0}$ and covariance $\boldsymbol{\Gamma}\!\left(\boldsymbol{\theta}\right)$ is
\begin{equation}\label{eq8_25_1}
p(\boldsymbol{y}; \boldsymbol{\theta}) = \frac{1}{(2\pi)^{\frac{N}{2}}\sqrt{\det \boldsymbol{\Gamma}\!\left(\boldsymbol{\theta}\right)}}\exp\left(-\frac{1}{2} \boldsymbol{y}^T \boldsymbol{\Gamma}\!\left(\boldsymbol{\theta}\right)^{-1}\boldsymbol{y}\right).
\end{equation}In MLE method, we seek the parameters $\boldsymbol{\theta}=(\alpha,\gamma, K, \ell)$ that would maximize the probability density function \eqref{eq8_25_1} with $\boldsymbol{y}$ being the observed velocity measures. Equivalently, we have to minimize the negative log of the likelihood function:
\begin{equation}\label{eq8_25_2}\begin{split}
NLL(\boldsymbol{\theta})&=-\log p(\boldsymbol{y}; \boldsymbol{\theta}) = \frac{1}{2}\boldsymbol{y}^T \boldsymbol{\Gamma}\!\left(\boldsymbol{\theta}\right)^{-1}\boldsymbol{y}+\frac{1}{2}\log \det \boldsymbol{\Gamma}\!\left(\boldsymbol{\theta}\right)+\frac{N}{2}\log (2\pi).\end{split}
\end{equation}Finding the minimum of the highly nonlinear function \eqref{eq8_25_2} with four parameters is computationally demanding. Therefore it is desirable to reduce the number of parameters which has to be estimated. Notice that the variance $s^2$ of  $\hat{Y}_{\alpha,\gamma}^{K, \ell}(t)$ is given by
$$s^2=\hat{C}_{\alpha,\gamma}^{K, \ell}(0)=\frac{K^2}{2\pi \alpha}\frac{\Gamma\left(\frac{1}{2\alpha}\right)
\Gamma\left(\gamma-\frac{1}{2\alpha}\right)}{\Gamma(\gamma)}.$$ As a result, the value of $K$ can be determined from this equation once $s^2, \alpha,\gamma$  are given. On the other hand, we can rewrite the covariance matrix $\boldsymbol{\Gamma}\!\left(\boldsymbol{\theta}\right)=\Gamma(\alpha,\gamma, K,\ell)$ as $s^2 \boldsymbol{\rho} (\boldsymbol{\theta}')= s^2\boldsymbol{\rho} (\alpha,\gamma,\ell)$, where $\boldsymbol{\rho} (\alpha,\gamma,\ell)$ is the correlation matrix $\left(\rho_{\alpha,\gamma}^{K,\ell}(i-j)\right)_{i, j=1}^{N}$, $$\rho_{\alpha,\gamma}^{K, \ell}(i-j) =\frac{C_{\alpha,\gamma}^{K,\ell}(i-j)}{C_{\alpha,\gamma}^{K,\ell}(0)},$$ which is independent of $K$. Rewriting in the variables $\alpha,\gamma, \ell, s^2$, we have
\begin{equation}\label{eq8_25_3}\begin{split}
NLL(\boldsymbol{\theta}',s^2)= &\frac{1}{2 s^2}\boldsymbol{y}^T \boldsymbol{\rho}\!\left(\boldsymbol{\theta}'\right)^{-1}\boldsymbol{y}+\frac{N}{2}\log s^2+\frac{1}{2}\log \det \boldsymbol{\rho}\!\left(\boldsymbol{\theta}\right)+\frac{N}{2}\log (2\pi).\end{split}
\end{equation}Taking derivative with respect to $s^2$, we find that for fixed $\boldsymbol{\theta}'=(\alpha,\gamma, \ell)$, the minimum of $NLL(\boldsymbol{\theta}',s^2)$ appears at \begin{equation}\label{eq8_25_5}s^2= \frac{1}{N}\boldsymbol{y}^T \boldsymbol{\rho}\!\left(\boldsymbol{\theta}'\right)^{-1}\boldsymbol{y}.\end{equation} Substituting this into \eqref{eq8_25_3}, we reduce the problem to finding $\boldsymbol{\theta}'$ to minimize the function
\begin{equation}\label{eq8_25_4}\begin{split}
\widetilde{NLL}(\boldsymbol{\theta}')= &\frac{N}{2}\log \boldsymbol{y}^T \boldsymbol{\rho}\!\left(\boldsymbol{\theta}'\right)^{-1}\boldsymbol{y}+\frac{1}{2}\log \det \boldsymbol{\rho}\!\left(\boldsymbol{\theta}'\right)+\frac{N}{2}\left(1+\log (2\pi)-\log(N)\right),\end{split}
\end{equation}and $s^2$ is then determined from \eqref{eq8_25_5}. The \emph{fminsearch}   function in Matlab which uses the simplex search algorithm by Neldon and Mead \cite{1_20_1}, is used to   identify the minimum of \eqref{eq8_25_4}. This algorithm does not involve computation of derivatives.  In order to compare the GWM model with that of WM, we also run the same search with $\alpha$ fix to $1$. The results are tabulated in Table \ref{Tab1}. It suggests the WM model $\hat{Y}_{1,\gamma}^{K,\ell}(t)$ with spectral density
\begin{equation*}
S_{\text{WM}}(\omega) = \frac{0.5796}{2\pi}\frac{1}{(\omega^2+0.75^2)^{1.02}}
\end{equation*}and the GWM model $\hat{Y}_{\alpha,\gamma}^{K,\ell}(t)$ with spectral density
\begin{equation*}
S_{\text{GWM}}(\omega) = \frac{50.9376}{2\pi}\frac{1}{(|\omega|^{1.03}+2.82^{1.03})^{4.12}}
\end{equation*}for the velocity measures.
From Table \ref{Tab1}, we see that the GWM model gives a better value to $\widetilde{NLL}(\boldsymbol{\theta}')$. On the other hand, a graphical comparison  of the PSD of the WM model and the GWM model for velocity measures and the empirical PSD (Figure \ref{Fig13}) also shows that the GWM model gives a better fit to the velocity measures compared to the WM model especially in the low frequency region.

\begin{table}
\caption{\label{Tab1} The estimated parameters for the WM and GWM models}

\begin{tabular}{||c|c|c|c|c|c|c||}
\hline
\hline
 &$\tilde{\alpha}$ & $\tilde{\gamma}$ & $\tilde{K}$ & $\tilde{\ell}$ &$\widetilde{s^2}$ &$\widetilde{NLL}$ \\
 \hline
 WM model & 1 & 1.0225 & 0.7857 & 0.7474 & 0.2994 & 1488.42\\
 GWM model & 0.5186 & 4.1223 &1.6965 & 2.8250 & 0.2995 & 1487.47\\
 \hline\hline
\end{tabular}
\end{table}

\begin{figure}
\epsfxsize=0.6\linewidth \epsffile{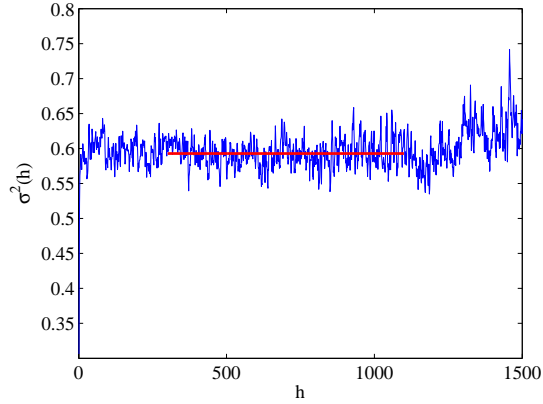} \caption{\label{Fig14}The empirical variogram $\tilde{\sigma}^2(h)$. The horizontal line gives the estimated value of $ \hat{C}_{\alpha,\gamma}^{K, \ell}(0)=0.2964$.}\end{figure}

\begin{figure}
\epsfxsize=0.6\linewidth \epsffile{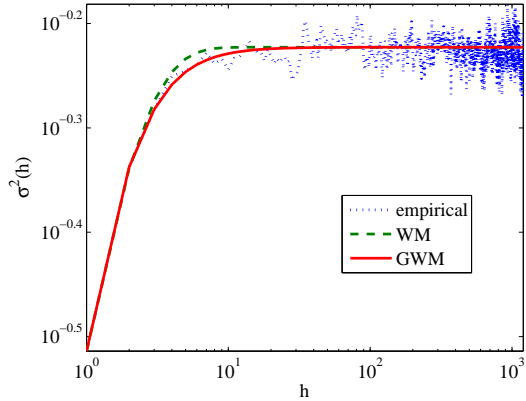} \caption{\label{Fig15}The empirical variogram $\tilde{\sigma}^2(h)$ compared to the variograms of the WM and GWM models.}\end{figure}

Here we would also like to remark that theoretically, the variogram $\hat{\sigma}_{\alpha,\gamma}^{K, \ell}(h)^2=\left\langle \left[\hat{Y}_{\alpha,\gamma}^{K,\ell}(t+h)-\hat{Y}_{\alpha,\gamma}^{K,\ell}(t)\right]^2\right\rangle=2(\hat{C}_{\alpha,\gamma}^{K, \ell}(0)-\hat{C}_{\alpha,\gamma}^{K, \ell}(h))$   approaches $$2\hat{C}_{\alpha,\gamma}^{K, \ell}(0)=\frac{K^2}{\pi \alpha}\frac{\Gamma\left(\frac{1}{2\alpha}\right)
\Gamma\left(\gamma-\frac{1}{2\alpha}\right)}{\Gamma(\gamma)}$$ as $h\rightarrow \infty$. Figure \ref{Fig14} shows the empirical variogram of the velocity measures estimated by $$\tilde{\sigma}^2 (h) =\frac{1}{N-h}\sum_{i=1}^{N-h} \left(y_{i+h}-y_i\right)^2.$$ The horizontal line gives an estimation of the variance $s^2 =0.2964$, which is very close to the one estimated by MLE. Figure \ref{Fig15} compares the empirical variogram to the variograms of the WM model and GWM model for velocity measures.

\section{ Concluding Remarks}
In this paper, we have introduced a new class of Gaussian random field with covariance belonging to a generalized Whittle-Mat$\acute{\text{e}}$rn family of covariance functions. Some of the basic properties of this GWM field are studied. Simulations of the GWM covariance and GWM field in two dimensions are carried out. We also apply this random process to model wind speed. It is shown that this new random field can provide a more flexible alternative to  modeling. In the future, we would like to extend the application of GWM field to provide models for other geostatistical data such as sea beam data, geothermal field temperature and soil data which the WM field has been shown to provide a good model \cite{10,11,12,13, 41}.

Just like its predecessor WM model,   GWM model will find its main applications in geostatistics. However, one expects it can have potential applications in modeling short range dependent process such as coding regions of DNA sequences and  fluctuations of an electropore of nano size \cite{le, lf, lg, lh, li, lj}. It will also be interesting to consider its applications in modeling fractional diffusion and fractional anomalous diffusion \cite{J5, nr12, J8, J9, J10, A1, A2}. This later aspect is being considered elsewhere.  For $d=2$, GWM may serve as a model to two-dimensional images arising in biology, chemistry and physics, in addition to those from geological and environmental images. The advances in imaging techniques allow better analysis of morphology of various material surfaces. Various spatial statistical and morphological methods are available to analyze the patterns of the surface of complex materials, hence the characterization of their physical properties. Examples of data that can be modeled by the random field include the concentration of particular component in a liquid or solid sample, properties such as porosity, permeability, conductivity, absorptivity,   emissivity, etc. of the material samples. Although applications of spatial models to statistical physics are still quite limited, some recent efforts have been made in \cite{1_5_1, 1_5_2, 1_5_3}. One expects GWM model also has such potential applications. In particular, as a correlation model, GWM model is useful in the modeling of morphological structure of complex material with spatial correlation that is short-ranged, that is, its underlying physical process is weakly correlated or weakly coupled over finite spatial or temporal scales. Readers can consult references in \cite{1_5_4} for recent advances and applications of spatial models in physics, in particular statistical physics and astrophysics.

	For modeling of data which may have   correlation with time or space dependent memory parameter or smoothness parameter, it is necessary to consider GWM field   $Y_{\alpha,\gamma}(\boldsymbol{t})$ with $\alpha$  and $\gamma$  replaced by $\alpha(\boldsymbol{t})>0$ and $\gamma(\boldsymbol{t})>0$, with $\alpha(\boldsymbol{t})\gamma\left(\boldsymbol{t}\right)>n/2$. A GWM field with two variable fractional indices can be studied in a similar way like the fractional Riesz-Bessel field with variable order \cite{28}. In fact, when $n=1$ and $\alpha=1$, such a generalization has been considered in \cite{n1} where it is called Weyl multifractional Ornstein-Uhlenbeck process. Another possible extension is the anisotropic counterpart of the GWM field whose covariance is a product of GWM processes.

\appendix
\section{Derivations of formulas (26) and (31)}
\noindent 1. We want to prove eq. \eqref{eq26}  when $\alpha\gamma\in
\left(\frac{n}{2}, \frac{n+2}{2}\right)$. Using regularization method, we have
\begin{equation*}\begin{split}
I=\lim_{a\rightarrow 0^+}& \Biggl\{\int_0^{\infty}
\frac{J_{\frac{n-2}{2}}(k)k^{\frac{n}{2}}}
{(k^2+a^2)^{\alpha\gamma}}dk-\frac{1}{2^{\frac{n-2}{2}}\Gamma\left(\frac{n}{2}\right)}
\int_0^{\infty}\frac{k^{n-1}dk}{(k^2+a^2)^{\alpha\gamma}}\Biggr\}.
\end{split}\end{equation*}Applying the formulas \#6.565, no. 4 and \#3.251, no. 11 of
\cite{34}, we find that
\begin{align*}
I=\lim_{a\rightarrow 0^+} \left\{
\frac{a^{\frac{n}{2}-\alpha\gamma}}{2^{\alpha\gamma-1}\Gamma(\alpha\gamma)}K_{\alpha\gamma-\frac{n}{2}}(a)
-\frac{a^{n-2\alpha\gamma}}{2^{\frac{n}{2}}}\frac{\Gamma\left(\alpha\gamma-\frac{n}{2}
\right)}{\Gamma(\alpha\gamma)} \right\}.
\end{align*}Now the formulas \#8.485 and \#8.445 of
\cite{34} give
\begin{equation}\label{eqA1}\begin{split}
&K_{\nu}(z)=K_{-\nu}(z) =\frac{\pi}{2\sin (\pi
\nu)}\left\{\sum_{j=0}^{\infty}
\frac{(z/2)^{2j-\nu}}{j!\Gamma(j+1-\nu)}-\sum_{j=0}^{\infty}\frac{(z/2)^{2j+\nu}}{j!
\Gamma(j+1+\nu)}\right\}
\end{split}\end{equation}when $\nu\in (0,1)$; which allows us to conclude that
\begin{equation*}\begin{split}
I=&-\frac{\pi}{2^{2\alpha\gamma-\frac{n}{2}}\sin\left[\pi\left(\alpha\gamma-\frac{n}{2}\right)\right]}
\frac{1}{\Gamma\left(\alpha\gamma\right)\Gamma\left(\alpha\gamma-\frac{n}{2}+1\right)}=
\frac{\Gamma\left(\frac{n}{2}-\alpha\gamma\right)}{2^{2\alpha\gamma-\frac{n}{2}}\Gamma(\alpha\gamma)}.\end{split}
\end{equation*}

\vspace{0.3cm}\noindent 2. We want to prove eq.~\eqref{eq31} when $\alpha\gamma=\frac{n+2}{2}$. Using \eqref{eq28}, we can write $I(\boldsymbol{t})$ as
the sum of $I_1(\boldsymbol{t})$ and $I_2(\boldsymbol{t})$, where
\begin{equation*}\begin{split}
I_1(\boldsymbol{t})=&\int_0^{1}\left(\frac{J_{\frac{n-2}{2}}(k)}
{k^{\frac{n-2}{2}}}-\frac{1}{2^{\frac{n-2}{2}}\Gamma\left(\frac{n}{2}\right)}
+\frac{k^2}{2^{\frac{n+2}{2}}\Gamma\left(\frac{n+2}{2}\right)}\right)
\frac{k^{n-1}}{(k^{2\alpha}+\lambda^2|\boldsymbol{t}|^{2\alpha})^{\gamma}}dk\\
&+\int_1^{\infty}\left(\frac{J_{\frac{n-2}{2}}(k)}
{k^{\frac{n-2}{2}}}-\frac{1}{2^{\frac{n-2}{2}}\Gamma\left(\frac{n}{2}\right)}\right)
\frac{k^{n-1}}{(k^{2\alpha}+\lambda^2|\boldsymbol{t}|^{2\alpha})^{\gamma}}dk\end{split}
\end{equation*}has a finite limit $I_1(0)$ as $|\boldsymbol{t}|\rightarrow 0$, and
\begin{equation*}
I_2(\boldsymbol{t}):=I(\boldsymbol{t})-I_1(\boldsymbol{t})=-\frac{1}{2^{\frac{n+2}{2}}\Gamma\left(\frac{n+2}{2}\right)}\int_0^1
\frac{k^{n+1}}{(k^{2\alpha}+\lambda^2|\boldsymbol{t}|^{2\alpha})^{\gamma}}dk.
\end{equation*} By making a change of variable $k\mapsto
k^{1/(2\alpha)}$, we have
\begin{equation*}
I_2(\boldsymbol{t})=-\frac{1}{2^{\frac{n+4}{2}}\alpha\Gamma\left(\frac{n+2}{2}\right)}\int_0^1\frac{
k^{\gamma-1}dk}{(k+\lambda^2|\boldsymbol{t}|^{2\alpha})^{\gamma}}.
\end{equation*} From this, we find that $I_2(\boldsymbol{t})$ can be written as a sum
of $I_3(\boldsymbol{t})$ and $I_4(\boldsymbol{t})$, where
\begin{equation*}\begin{split}
I_3(\boldsymbol{t})
=&-\frac{1}{2^{\frac{n+4}{2}}\alpha\Gamma\left(\frac{n+2}{2}\right)}  \int_0^1\left\{\frac{
k^{\gamma-1}}{(k+\lambda^2|\boldsymbol{t}|^{2\alpha})^{\gamma}}-\frac{1}{k+\lambda^2|\boldsymbol{t}|^{2\alpha}}
\right\}dk\end{split}
\end{equation*}has a finite limit $I_3(\mathbf{0})$ when $|\boldsymbol{t}|\rightarrow 0$, and
 \begin{equation*}\begin{split}
I_4(\boldsymbol{t}):=&I_2(\boldsymbol{t})-I_3(\boldsymbol{t})=-\frac{1}{2^{\frac{n+4}{2}}\alpha\Gamma\left(\frac{n+2}{2}\right)}\int_0^1
\frac{1}{k+\lambda^2|\boldsymbol{t}|^{2\alpha}}dk\\=&-\frac{1}{2^{\frac{n+4}{2}}\alpha\Gamma\left(\frac{n+2}{2}\right)}\log
\frac{1+\lambda^2|\boldsymbol{t}|^{2\alpha}}{\lambda^2|\boldsymbol{t}|^{2\alpha}}\\
=&-\frac{1}{2^{\frac{n+2}{2}}\Gamma\left(\frac{n+2}{2}\right)}\log\frac{1}{|\boldsymbol{t}|}
+\frac{1}{2^{\frac{n+4}{2}}\alpha\Gamma\left(\frac{n+2}{2}\right)}\log\lambda^2+o(1).\end{split}
\end{equation*}  Therefore,
we have shown that
$$I(\boldsymbol{t})
=-\frac{1}{2^{\frac{n+2}{2}}\Gamma\left(\frac{n+2}{2}\right)}\log\frac{1}{|\boldsymbol{t}|}+A+o(1),$$where
$$A=I_1(\mathbf{0})+I_3(\mathbf{0})+\frac{1}{2^{\frac{n+4}{2}}\alpha\Gamma\left(\frac{n+2}{2}\right)}\log\lambda^2.$$

\vspace{1cm}\noindent\textbf{Aknowledgments}\;
The authors would
like to thank Malaysian Academy of Sciences, Ministry of Science,
Technology  and Innovation for funding this project under the
Scientific Advancement Fund Allocation (SAGA) Ref. No P96c.

\end{document}